%% file: Frittelli_et_al_arXiv.tex
\documentclass[a4paper, 11pt]{article}
\usepackage[T1]{fontenc}
\usepackage[utf8]{inputenc}
\usepackage[english]{babel}
\usepackage{amsmath}
\usepackage{amssymb}
\usepackage{amsthm}
\usepackage{color}
\usepackage{graphicx}
\usepackage{tikz}
\usepackage{placeins}
\usepackage{epstopdf}
\usepackage{subfigure}
\usepackage{mathtools}
\usepackage{caption}
\usepackage{enumitem}
\usepackage{rotating}
\usepackage{hyperref}
\usepackage{authblk}
\usepackage{soul}
\usepackage[margin = 2.8cm]{geometry}

\newtheorem{theorem}{Theorem}
\newtheorem{lemma}{Lemma}
\newtheorem{definition}{Definition}

\newcommand{\red}[1]{\textcolor{red}{#1}}

\newcommand{\bolda}{\mathbf{a}}

\newcommand{\boldf}{\mathbf{f}}

\newcommand{\boldu}{\mathbf{u}}
\newcommand{\boldU}{\mathbf{U}}
\newcommand{\boldv}{\mathbf{v}}
\newcommand{\boldV}{\mathbf{V}}
\newcommand{\boldx}{\mathbf{x}}
\newcommand{\boldy}{\mathbf{y}}
\newcommand{\boldxi}{\boldsymbol{\xi}}

\newcommand{\boldphi}{\boldsymbol{\phi}}
\newcommand{\boldvarphi}{\boldsymbol{\varphi}}

\newcommand{\boldtheta}{\boldsymbol{\theta}}
\newcommand{\boldnu}{\boldsymbol{\nu}}
\newcommand{\boldrho}{\boldsymbol{\rho}}

\newcommand{\norm}[1]{\left\|{#1}\right\|}
\newcommand{\semin}[1]{\left|{#1}\right|}

\newcommand{\intgamma}{\int_{\Gamma}}
\newcommand{\intgammah}{\int_{\Gamma_h}}
\newcommand{\nablagamma}{\nabla_\Gamma}
\newcommand{\nablagammah}{\nabla_{\Gamma_h}}

\title{Lumped finite element method for reaction-diffusion systems on compact surfaces}
\author[1]{Massimo Frittelli}
\author[2]{Anotida Madzvamuse}
\author[1]{Ivonne Sgura}
\author[3]{Chandrasekhar Venkataraman}

\affil[1]{\small{Dipartimento di Matematica e Fisica “E. De Giorgi”, Università del Salento, via per Arnesano, I-73100 Lecce, Italy}}
\affil[2]{\small{University of Sussex, School of Mathematical and Physical Sciences, Department of Mathematics, University of Sussex, Brighton, BN1 9QH, United Kingdom}}
\affil[3]{\small{School of Mathematics and Statistics, University of St Andrews, Fife, KY16 9SS, United Kingdom}}

\date{}

\numberwithin{equation}{section}
\begin{document}

\maketitle

\begin{abstract}
We propose and analyse a novel surface finite element method that preserves the invariant regions of systems of semilinear parabolic equations on closed compact surfaces  in $\mathbb{R}^3$ under discretisation. We also provide a fully-discrete scheme by applying the implicit-explicit
(IMEX) Euler method in time. We prove the preservation of the invariant
rectangles of the continuous problem under spatial and full discretizations.
For scalar equations, these results reduce to the well-known discrete maximum principle. 
Furthermore, we prove optimal error bounds for the semi- and fully-discrete methods, that is the convergence rates are quadratic in the meshsize and
linear in the timestep.
Numerical experiments are provided to support the theoretical findings.
In particular we provide examples in which, in the absence of lumping, the numerical solution violates the invariant region leading to blow-up due to the nature of the kinetics.
\end{abstract}

\section{Introduction}
\label{intro}

Partial differential equations (PDEs) of the form of reaction-diffusion systems (RDSs) have been extensively employed to model many different processes in a wide range of fields such as biology \cite{murray2001, kondo1995reaction, nijhout2003pigmentation, ferreira2002reaction}, chemistry \cite{castets1990, vanag2004waves}, electrochemistry \cite{bessler2005new, lacitignola2015} and finance \cite{mainardi2000fractional, becherer2005}. In many applications the domain of integration is a stationary or an evolving curved surface, rather than a planar region. For instance, surface RDSs have been applied to the study of biological patterning \cite{barreira2011surface}, tumour growth \cite{chaplain2001}, superconductivity \cite{du2005approximations}, metal dealloying \cite{eilks2008numerical}, biomembrane modeling \cite{elliott2010modeling}, cell motility \cite{elliott2012modelling} and phase separation \cite{tang2005phase}, just to mention a few examples. In this paper we consider RDSs of arbitrarily many equations on a stationary surface of the form:
 for $i=1,\dots,r,$
\begin{equation}
\label{continuoussystem}
\begin{cases}
&\dfrac{\partial u_i}{\partial t} - d_i \Delta_\Gamma u_i = f_i(u_1,\dots, u_r), \qquad\text{ in }\Gamma\times(0,T],\\
&u_i(\boldx,0) = {u_0}_i (\boldx),\qquad \forall\ \boldx  \in\Gamma,
\end{cases}
\end{equation}
where $\Gamma$ is a smooth stationary orientable surface of codimension one in $\mathbb{R}^3$ without boundary, $\Delta_\Gamma$ is the Laplace-Beltrami operator on $\Gamma$  (which is defined as the tangential divergence of the tangential gradient, see \cite{acta} for the definitions), $d_i$ are strictly positive diffusion coefficients and ${u_0}_i$ are smooth, bounded functions.

A key feature of many RDSs is the existence of invariant regions. A region $\Sigma$ in the phase space $\mathbb{R}^r$ is said to be invariant for \eqref{continuoussystem} if, whenever the initial condition has values in $\Sigma$, the solution of \eqref{continuoussystem} stays in $\Sigma$ as long as it is defined. Knowing that a given model possesses an invariant region is useful in a couple of ways. First, when solving RDSs arising from applications, solutions are usually meaningful as long as they range within a limited set of values (an example is given by mass-action laws \cite{chellaboina2009modeling}, in which solutions are required to be componentwise nonnegative). Second, an invariant region provides an a-priori bound on the analytical solution which can be helpful, for instance, when studying the convergence of numerical methods. Sufficient conditions for a region to be invariant for a given RDS were given in \cite{smollerbook,smollerarticle} on planar domains and in \cite[p. 335-353]{taylor1997partial} on stationary surfaces. In both cases, for distinct $d_i$'s, the only possible invariant regions for \eqref{continuoussystem} are (bounded or unbounded) hyper-rectangles in $\mathbb{R}^r$, that is to say regions in the form
\begin{equation}
\label{invariantrectangle}
\Sigma = \prod_{i=1}^r [m_i, M_i],
\end{equation}
with $m_i,M_i\in\mathbb{R}\cup \{-\infty, +\infty\}$ for all $i=1,\dots,r$, whereas if some $d_i$ coincide, more general regions are allowed to be invariant. Since we are addressing general diffusion coefficients, we will consider invariant hyper-rectangles 
\eqref{invariantrectangle} only. Among the literature RD models having invariant hyper-rectangles we recall the Gierer-Meinhardt \cite{kovacs2003spatial}, Hodgkin-Huxley \cite{galusinski1998}, FitzHugh-Nagumo \cite{rauch1978}, Oregonator \cite{you2011}, Rosenzweig-Macarthur \cite{skalski2001functional, gonzalez2003dynamic}, and the spatially extended Lotka-Volterra \cite{alikakos1979} models.
For $r=1$ in \eqref{continuoussystem}, i.e. scalar semilinear equations, we remark that the min-max condition and the maximum principle, given by
\begin{align}
\label{minmaxcondition}
\min_{\Gamma} u_0 \leq u(\boldx,t) \leq \max_\Gamma u_0,\qquad \forall\ (\boldx,t)\in \Gamma\times [0,T],\\
0 \leq u(\boldx,t) \leq \max_\Gamma u_0,\qquad \forall\ (\boldx,t)\in \Gamma\times [0,T],
\end{align}
respectively, correspond to particular invariant regions, given by\\
$\Sigma = [\min_{\Gamma} u_0, \max_\Gamma u_0]$ and $\Sigma = [0, \max_\Gamma u_0]$, respectively.

The increasing interest from applications in RDSs on manifolds has stimulated the development of numerical methods for such systems. Among the methods for PDEs on stationary surfaces we recall: finite differences \cite{varea1999turing}, the spectral method of lines \cite{chaplain2001}, closest point methods (see \cite{macdonald2009implicit} and references therein), kernel methods (see \cite{fuselier2015radial} and references therein), embedding methods (see \cite{bertalmio2003variational} and references therein), and surface finite element methods (SFEM) (see \cite{dziuk1988, acta, tuncer2015} and references therein). In this paper we consider a lumped mass surface finite element method (LSFEM) for the spatial discretization of Eqs. \eqref{continuoussystem}. We recall that finite elements with mass lumping have already been applied to reaction-diffusion systems on planar domains, see for example \cite{thomee1980, garvie2007}.

 To carry out a \emph{fully-discrete scheme} we will follow an implicit-explicit (IMEX) approach, i.e. by treating diffusion implicitly and reactions explicitly. Among the class of IMEX methods, we will consider the simplest one, the IMEX Euler scheme considered for example in \cite{madzvamuse2006time, lakkis2013implicit}. IMEX methods have been widely applied in fluid dynamics, combined with spectral methods on planar domains  \cite{canuto1988spectral, kim1985application}, in reaction-diffusion problems, in combination with finite differences in space on planar domains \cite{ruuth1995implicit}, with finite elements on stationary planar domains \cite{elliott1993global}, on evolving planar domains \cite{madzvamuse2006time}, and with the closest point method on stationary surfaces \cite{macdonald2008closest}. An error analysis of finite element approximations with IMEX timestepping for semilinear systems on evolving domains is carried out in \cite{lakkis2013implicit}.

When numerically approximating RDSs, it is desirable for the considered numerical method to preserve invariant rectangles of the continuous problem. 
For the scalar case, that is the maximum principle, on planar domains, works in this direction cover the homogeneous heat equation (see \cite{newthomee} and references therein), reaction-diffusion equations \cite{thomee1980, farago2012discrete, elliott1993global}, anisotropic reaction-diffusion \cite{li2013maximum} and reaction convection-diffusion equations \cite{lu2014maximum}. For reaction-diffusion systems of many equations on planar domains, the problem is addressed in \cite{hoff1978stability}. 
The aforementioned works consider different spatial methods. Most of them require the disctretisation to be sufficiently refined, in order to preserve invariant rectangles and maximum principles. A notable exception is the lumped finite element method (LFEM) \cite{elliott1993global, thomee1980, newthomee, li2013maximum, lu2014maximum}.
To the best of our knowledge, numerical methods for surface RDSs preserving the invariant rectangles of the continuous problem have not yet been presented. As far as we know, only a time dependent discrete maximum principle for a scalar diffusion problem on evolving surfaces is given in a recent work \cite{kovacs2015}, in which the evolving surface finite element method (ESFEM) is applied. This motivates the present study in which we introduce the LSFEM, which not only preserves the invariant rectangles at the discrete level, but also requires no restriction on the mesh size.

The main contributions of this paper are twofold. First, we prove discrete maximum principles for the LSFEM semi-discretisation and IMEX Euler-LSFEM full discretisation for a class of semilinear parabolic equations, and the preservation of invariant rectangles under discretisation for weakly coupled (i.e. coupled only through the reaction kinetics) semilinear (i.e. in which only the kinetics are nonlinear) RDSs \eqref{continuoussystem}. Second, we prove optimal error bounds for the semi-discrete and fully discrete schemes. Among the numerical tests, we provide an example of RDS possessing an invariant region, in which the SFEM blows-up, while the LSFEM preserves the region.

The present article is structured as follows. In Section \ref{sec:heatproblem} we consider a semilinear scalar parabolic equation on a closed orientable surface in strong and weak formulation, we present its LSFEM space discretization, its Euler IMEX/LSFEM full discretization and prove the preservation of the maximum principle under spatial and full discretization in Theorems \ref{thm:sdmaximumprinciple} and \ref{thm:fdmaximumprinciple}, respectively. In Section \ref{sec:rdsysyems} we consider a general RDS of arbitrarily many equations on closed orientable surfaces, we derive its LSFEM space discretization, its Euler IMEX/LSFEM time discretization and prove the preservation of the invariant rectangles under spatial and full discretizations in Theorems \ref{thm:sdinvariantregion} and \ref{thm:fdinvariantregion}, respectively. In Section \ref{sec:erroranalysis}, optimal error estimates for both the semi- and fully-discrete methods introduced in the previous sections are proven in Theorems \ref{thm:SDerroranalysis} and \ref{thm: FDerroranalysis}, respectively. Numerical experiments are shown in Section \ref{sec:numericaltest}.

\section{A semilinear scalar parabolic equation}
\label{sec:heatproblem}
\subsection{The continuous problem}
\label{sec:continuousproblem}

We start by considering scalar parabolic PDEs in order to illustrate the main ideas behind the approach described in this work and to introduce the analysis in a less technical setting.

Let $\Gamma$ be a compact, orientable, smooth surface of codimension one in $\mathbb{R}^3$ without boundary. We assume that $\Gamma$ is represented as the zero level set of a sufficiently smooth \emph{signed distance function} $d$, defined in an open neighbourhood $W$ of $\Gamma$ such that $\nabla d(\boldx) \neq \boldsymbol{0}\ \forall \boldx\in W$ by
\begin{equation*}
 \Gamma = \{\boldx\in W| d(\boldx)=0\}.
\end{equation*}
The normal unit vector on $\Gamma$ is then defined by
\begin{equation*}
\boldnu(\boldx) = \frac{\nabla d(\boldx)}{|\nabla d(\boldx)|},\qquad \forall\boldx\in\Gamma.
\end{equation*}
We assume that every point $\boldx\in W$ may be uniquely represented as
\begin{equation}
\label{onetoone}
\boldx = \bolda(\boldx) + d(\boldx)\boldnu(\bolda(\boldx)),
\end{equation}
with $\bolda(\boldx)\in\Gamma$. A sufficient condition on the thickness of $W$ such that this property holds is given in \cite{acta}.

We briefly recall the definitions of Sobolev and Bochner spaces on surfaces. For $q\in \mathbb{N}\cup \{0\}$, the Sobolev space $H^q(\Gamma)$ is the space of functions $u:\Gamma\rightarrow\mathbb{R}$ such that, for all $i=0,\dots,q$, the $i$-th order tangential derivatives, meant in a distributional sense, are $L^2(\Gamma)$, whilst $H^{-q}(\Gamma)$ is the dual space of $H^q(\Gamma)$, that is the space of linear continuous functionals on $H^q(\Gamma)$.  For $p\in [1,+\infty]$, if $X$ is a Banach space, the Bochner space $L^p([0,T]; X)$ is the space of functions $u:[0,T]\rightarrow X$ such that the function $\|u\|_{X}: [0,T]\rightarrow\mathbb{R}$ is $L^p([0,T])$. For further details on Sobolev and Bochner spaces on surfaces we refer the interested reader to \cite{gilbarg2015elliptic}, \cite{hebey2008sobolev} or \cite{taylor1997partial}.

In this section we consider the following semilinear parabolic equation posed on $\Gamma$:
\begin{equation}
\label{systemunderstudy}
\dot{u} - d\Delta_\Gamma u = - \beta u^\alpha, \qquad \boldx\in\Gamma,\ t\in (0,T],
\end{equation}
where the dot denotes the time derivative, $d>0$, $\alpha \geq 1$, $\beta\geq 0$, endowed with the nonnegative $\mathcal{C}^2(\Gamma)$ initial condition
\begin{equation*}
u(\boldx,0) = u_0(\boldx),\qquad \boldx\in\Gamma.
\end{equation*}
The requirement that $\alpha\geq 1$ is needed to make the source term $u^\alpha$ be Lipschitz in a neighbourhood of $u=0$, which is a necessary condition for the existence and uniqueness of a solution at all positive times. The conditions $\beta\geq 0$ and $u_0\geq 0$ together are needed to guarantee the maximum principle \eqref{minmaxcondition}. The homogeneous heat equation is obtained as a special case for $\beta=0$. The weak formulation of the problem seeks to find a $u\in L^2([0,T]; H^1(\Gamma)) \cap L^\infty([0,T]\times \Gamma)$ with $\dot{u}\in L^2([0,T]; H^{-1}(\Gamma))$ such that
\begin{equation}
\label{weakformulation}
\int_{\Gamma} \dot{u}\varphi + d\int_{\Gamma}\nabla_\Gamma u\cdot\nabla_\Gamma \varphi = -\beta\int_\Gamma u^\alpha\varphi,\qquad \forall\ \varphi\in H^1(\Gamma).
\end{equation}

\subsection{Space discretization}
\label{sec:spacediscterization}
As mentioned previously, in the present work our focus is on finite element discretisations. We now present the necessary notation and concepts needed to describe the numerical method.

Given $h>0$, a triangulated surface $\Gamma_h\subset W$ is defined by
\begin{equation*}
\Gamma_h = \bigcup_{K\in\mathcal{K}_h} K,
\end{equation*}
where $\mathcal{K}_h$ is a set of finitely many non degenerate triangles, whose diameters do not exceed $h$ and whose vertices $\{\boldx_i\}_{i=1}^N$ lie on $\Gamma$, such that, for $\bolda(\boldx)$ as defined in \eqref{onetoone}, $\bolda_{|\Gamma_h}(\boldx)$ is a one-to-one map between $\Gamma$ and $\Gamma_h\subset W$.

Following \cite{acta}, we define lifts and unlifts. Given a function $V: \Gamma_h\rightarrow\mathbb{R}$, its \emph{lift} $V^\ell:\Gamma\rightarrow\mathbb{R}$ is defined by
\begin{equation*}
V^\ell(\bolda(\boldx)) = V(\boldx),\qquad \forall \boldx\in \Gamma_h.
\end{equation*}
Given a function $v:\Gamma\rightarrow\mathbb{R}$, its \emph{unlift} $v^{-\ell}:\Gamma_h\rightarrow\mathbb{R}$ is defined by
\begin{equation*}
v^{-\ell}(\boldx) = v(\bolda(\boldx)),\qquad \forall\ \boldx\in \Gamma_h.
\end{equation*}
Next, let $S_h$ be the space of piecewise linear functions on $\Gamma_h$ defined by
\begin{equation*}
S_h = \{V\in\mathcal{C}^0(\Gamma_h)\ |\ {V}_{|K} \text{ is linear affine } \forall K\in\mathcal{K}_h\}
\end{equation*}
and $S_h^\ell$ be its lifted counterpart:
\begin{equation*}
S_h^\ell = \{V^\ell\ |\ V\in S_h \}.
\end{equation*}
Let $\{\chi_i\}_{i=1}^N$ be the nodal basis of $S_h$ defined by
\begin{equation}
\label{nodalbasis}
\chi_i(\boldx_j) = \delta_{ij},\qquad \forall i,j=1,\dots,N.
\end{equation}
For $v\in\mathcal{C}^0(\Gamma_h)$, the piecewise linear interpolant $I_h(v)$ of $v$ is the function in $S_h$ given by
\begin{equation}
\label{linearinterpolant}
I_h(v) = \sum_{i=1}^N v(\boldx_i)\chi_i.
\end{equation}
\noindent
We define the following semi-discrete problem: find $U\in L^2([0,T]; S_h)$ with $\dot{U}\in L^2([0,T]; S_h)$ such that
\begin{equation}
\label{semi-discreteformulation}
\intgammah I_h(\dot{U}\phi) +d\intgammah\nablagammah U\cdot\nablagammah\phi = -\beta\intgammah I_h(U^\alpha\phi),\qquad \forall\ \phi\in S_h.
\end{equation}
We express the semi-discrete solution $U$ in terms of the nodal basis \eqref{nodalbasis} as
\begin{equation}
\label{nodalbasisexpansion}
U(\boldx,t) = \sum_{i=1}^N \xi_{i}(t)\chi_i(\boldx). 
\end{equation}
We then define the lumped mass matrix $\bar{M}=(\bar{m}_{ij})$ and the stiffness matrix $A = (a_{ij})$, respectively, by
\begin{align}
\label{lumpedmassmatrix}
&\bar{m}_{ij} = \int_{\Gamma_h} I_h(\chi_i\chi_j),\qquad \forall\ i,j=1,\dots, N,\\
\label{stiffnessmatrix}
&a_{ij} = \int_{\Gamma_h} \nabla_{\Gamma_h}\chi_i\cdot\nabla_{\Gamma_h}\chi_j,\qquad\forall \ i,j=1,\dots,N.
\end{align}
We recall that the mass matrix used in the standard SFEM \cite{dziuk1988, acta} is defined by
\begin{equation*}
m_{ij} = \int_{\Gamma_h} \chi_i\chi_j,\qquad \forall\ i,j=1,\dots, N.\\
\end{equation*}
Hence, the semi-discrete problem \eqref{semi-discreteformulation} can be expressed as the following ODE system:
\begin{equation}
\label{sys1}
\bar{M}\dot{\boldxi} + dA\boldxi = -\beta\bar{M}\boldxi^\alpha,
\end{equation}
where $\boldxi := (\xi_1,\dots, \xi_N)^T$. In the following we will show that, under suitable assumptions on the triangulation $\mathcal{K}_h$, this method fulfills a discrete maximum principle, that is the discrete version of \eqref{minmaxcondition}. To this end we introduce a regularity assumption for the triangulation on the mesh $\mathcal{K}_h$ which mimicks the standard Delaunay condition on planar domains and then we show how it affects the properties of the stiffness matrix $A$ in \eqref{stiffnessmatrix}.

\begin{figure}
\centering
\input{figure1.tex}
\caption{Schematic representation of condition \eqref{angle_condition} for triangles $K_1$ and $K_2$.}\label{fig:diagram}
\end{figure}
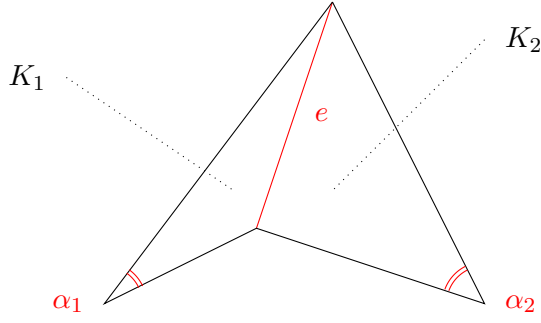

Let $e$ be an edge of the triangulation $\mathcal{K}_h$ and let $K_1$ and $K_2$ be the triangles sharing the edge $e$. Let $\alpha_1$ and $\alpha_2$ be the angles in $K_1$ and $K_2$ opposite to $e$, respectively. For every edge $e$ in $\mathcal{K}_h$ we require that
\begin{equation}
\label{angle_condition}
\alpha_1 + \alpha_2 \leq \pi.
\end{equation}
This condition is represented in Fig. \ref{fig:diagram}.
The following result extends to triangulated surfaces a characterization of \eqref{angle_condition} given in \cite{thomeebook} for the planar case.
\begin{lemma}
\label{delaunaylemma}
$\mathcal{K}_h$ fulfills \eqref{angle_condition} if and only if
\begin{equation}
\label{delaunayproperty}
(\nabla_{\Gamma_h}\chi_i,\nabla_{\Gamma_h}\chi_j) \leq 0\qquad \forall\ i\neq j.
\end{equation}
\begin{proof}
Let $\boldx_i$ and $\boldx_j$ be two distinct nodes of $\mathcal{K}_h$. If $\boldx_i$ and $\boldx_j$ are not neighbours, then $(\nabla_{\Gamma_h}\chi_i,\nabla_{\Gamma_h}\chi_j) = 0$. Otherwise, let $e$ be the edge connecting $\boldx_i$ and $\boldx_j$. Since the intersection of the support of the pyramidal functions $\chi_i$ and $\chi_j$ is $K_1 \cup K_2$ (see Fig. \ref{fig:diagram}) then we can write
\begin{equation}
\label{delaunayexpression}
(\nabla_{\Gamma_h}\chi_i,\nabla_{\Gamma_h}\chi_j) = (\nabla_{K_1}\chi_i,\nabla_{K_1}\chi_j) + (\nabla_{K_2}\chi_i,\nabla_{K_2}\chi_j).
\end{equation}
Let $T_1$ and $T_2$ be two direct isometries (that is with $\det(J_1)=\det(J_2)=1$) that map $K_1$ and $K_2$ into two triangles $K_1^0$ and $K_2^0$ contained in the $xy$ plane, respectively, and let $J_1$ and $J_2$ be the Jacobians of $T_1$ and $T_2$, respectively. Then, expression \eqref{delaunayexpression} can be written equivalently as 
\begin{equation}
\label{delaunayexpression2}
\begin{split}
&\int_{K_1^0} \left(J_1\nabla_{K_1^0}(\chi_i\circ T_1^{-1}))\cdot (J_1\nabla_{K_1^0}(\chi_j\circ T_1^{-1})\right)\det(J_1)+\\
&\int_{K_2^0} \left(J_2\nabla_{K_2^0}(\chi_i\circ T_2^{-1}))\cdot (J_2\nabla_{K_2^0}(\chi_j\circ T_2^{-1})\right)\det(J_2).
\end{split}
\end{equation}
Since $\det(J_1)=\det(J_2)=1$ and $\nabla_{K_1^0}$ and $\nabla_{K_2^0}$ both collapse to the standard gradient $\nabla$ in $\mathbb{R}^2$, expression \eqref{delaunayexpression2} then becomes
\begin{equation*}
\int_{K_1^0} \nabla(\chi_i\circ T_1^{-1})\cdot \nabla(\chi_j\circ T_1^{-1})+\int_{K_2^0} \nabla(\chi_i\circ T_2^{-1})\cdot \nabla(\chi_j\circ T_2^{-1}).
\end{equation*}
It is known that (see \cite{thomeebook}) this expression only depends on the transformed angles $\alpha_1^0 = \alpha_1$, $\alpha_2^0=\alpha_2$ and is given by
\begin{equation}
-\frac{\sin(\alpha_1+\alpha_2)}{4\sin(\alpha_1)\sin(\alpha_2)},
\end{equation}
which is nonpositive if and only if $\alpha_1+\alpha_2\leq\pi$. This completes the proof.
\end{proof}
\end{lemma}

Now, let $\vec{1}$ and $\vec{0}$ be the vector of ones and the null vector in $\mathbb{R}^N$, respectively. As shown in \cite{thomeebook} (pages 272-273), the structure \eqref{delaunayproperty} of the stiffness matrix, together with the diagonal structure \eqref{lumpedmassmatrix} of the lumped mass matrix, imply that, for every $s>0$, $\bar{M}+sA$ is an M-matrix. It then follows that
\begin{equation}
\label{matrixproperty1}
(\bar{M}+sA)^{-1}\bar{M} \geq \vec{0},
\end{equation}
meaning that this matrix has nonnegative entries. If $\boldxi = \vec{1}$, from \eqref{nodalbasisexpansion} we have $U(\boldx,t)=1$ for all $(\boldx,t)\in\Gamma_h\times [0,T]$, and thus $\nablagammah U(\boldx,t)$ vanishes, which yields $A\vec{1}=\vec{0}$. It therefore follows that
\begin{equation}
\label{matrixproperty2}
(\bar{M}+sA)^{-1}\bar{M}\vec{1} = \vec{1}.
\end{equation}
We will show that \eqref{matrixproperty1} and \eqref{matrixproperty2} play a crucial role in the discrete maximum principle for the parabolic equation \eqref{systemunderstudy} and the preservation of invariant regions of reaction-diffusion systems (see next Section \ref{sec:rdsysyems}).

\subsection{Time discretization}
\label{sec:timediscretization}
By applying the IMEX Euler scheme (i.e. treating diffusion implicitly and the reactions explicitly), with time step $\tau>0$, to \eqref{sys1} we obtain the fully-discrete scheme
\begin{equation}
\label{imexeuler}
\displaystyle\bar{M}\frac{\boldxi^{n+1}-\boldxi^n}{\tau} + dA\boldxi^{n+1} = -\beta\bar{M}(\boldxi^n)^\alpha,\qquad \forall\ n\in\mathbb{N}\cup \{0\},
\end{equation}
with $\boldxi^0 = \boldxi(0)$, where $\boldxi(t)$ is defined in \eqref{nodalbasisexpansion}, or equivalently,
\begin{equation}
\label{fullydiscretescheme}
\boldxi^{n+1} = (\bar{M}+d\tau A)^{-1}\bar{M}(\boldxi^n -\tau\beta(\boldxi^n)^\alpha), \qquad \forall\ n\in\mathbb{N}\cup \{0\}.
\end{equation}
We remark that, for $\beta=0$ (for the case of the homogeneous heat equation), the timestepping scheme collapses to the standard implicit Euler method.

\subsection{Semi- and fully-discrete maximum principles}
\label{sec:maximumprinciple}
It is known that the lumped FEM fulfills a discrete maximum principle for the homogeneous heat equation on planar domains, see \cite{thomeebook}. This result has been generalized to general diffusion problems in divergence form in \cite{thomee1980}. The purpose of this section is to extend this result to equation \eqref{systemunderstudy}, which includes, as a special case, the homogeneous heat equation on $\Gamma$.
\begin{theorem}[Maximum principle for \eqref{semi-discreteformulation}]
\label{thm:sdmaximumprinciple}
The semi-discrete solution $\boldxi(t)$ of \eqref{semi-discreteformulation} fulfills the following maximum principle
\begin{equation}
0 \leq \xi_i(t) \leq \max_{\mathbb{R}^N} \left\lbrace \boldxi(0)\right\rbrace, \qquad \forall\ i=1,\dots,N,\ \forall\ t>0.
\end{equation}
\begin{proof}
We rewrite \eqref{sys1} as
\begin{equation}
\label{sys2}
\dot{\boldxi} + d\bar{M}^{-1}A\boldxi = -\beta\boldxi^\alpha.
\end{equation}
Consider the auxiliary equation
\begin{equation}
\label{auxiliarysystem}
\dot{\boldxi} = -d\bar{M}^{-1}A\boldxi -\beta|\boldxi|^\alpha\emph{sign}(\boldxi),
\end{equation}
where $|\boldxi|$ and $\emph{sign}(\boldxi)$ are the componentwise absolute value and the componentwise sign function of $\boldxi$, respectively. If $\mu = \max_{\boldx\in\Gamma_h}v_h(\boldx)$, it is sufficient to prove that the solution of the ODE system \eqref{sys2} does not escape the region  $\Sigma = [0,\mu]^N$, i.e. we have to prove that, for every $\varepsilon > 0$, the solution of \eqref{auxiliarysystem} does not leave the region $\bar{\Sigma} := [-\varepsilon, \mu]^N$. To this end, we have to prove that the vector field on the right-hand-side of \eqref{auxiliarysystem}, computed on every $(N-1)$-dimensional face of $\bar{\Sigma}$, points toward the interior of $\bar{\Sigma}$. To this end, let $\boldxi$ be a point on $\partial \bar{\Sigma}$. This means that there exists $i=1,\dots,N$ such that $\xi_{i}\in \{-\varepsilon,\mu\}$. Suppose $\xi_{i} = \mu$; in the case $\xi_{i} = -\varepsilon$ the proof is analogous. Then
\begin{equation}
\label{boundaryheat}
\xi_{j}\leq \xi_{i},\quad j\neq i.
\end{equation}
All we have to prove is that $\dot{\xi}_{i}$ is negative. To this end, we prove that:
\begin{enumerate}
\item $-|\xi_i|^\alpha\text{sign}(\xi_i) = -|\mu|^\alpha \text{sign}(\mu) < 0$ from \eqref{boundaryheat};
\item The $i^{th}$ component of the vector $-d\bar{M}^{-1}A\boldxi$ is nonpositive. In fact, since $\bar{M}$ is a diagonal matrix, this component is given by
\begin{equation}
\label{expression1heat}
-(d\bar{M}^{-1}A\boldxi)_i = -d\bar{m}_{ii}^{-1}\sum_{j=1}^N a_{ij}\xi_{j}.
\end{equation}
We can split the sum on the right-hand-side by isolating the $a_{ii}\xi_{i}$ term:
\begin{equation}
\label{expression2heat}
d\bar{m}_{ii}^{-1}\left(-a_{ii}\xi_{i} +\sum_{j\in \{1,\dots,N\}\setminus \{i\}} (-a_{ij})\xi_{j}\right).
\end{equation}
Since $a_{ij}\leq 0$ for $i\neq j$ from Lemma \ref{delaunaylemma} and $\xi_{j}\leq \xi_{i}$ for $j\neq i$ from \eqref{boundaryheat}, expression \eqref{expression2heat} is less than or equal
\begin{equation}
\label{expression3heat}
\begin{split}
&d\bar{m}_{ii}^{-1}\xi_{i}\left(-a_{ii} +\sum_{j\in\{1,\dots,N\} \setminus \{ i\}} (-a_{ij})\right) =-d\bar{m}_{ii}^{-1}\xi_{i}\sum_{j=1}^N a_{ij}.
\end{split}
\end{equation}
From the definition of $A$, the right hand side of \eqref{expression3heat} is equal to
\begin{equation}
\label{expression4heat}
d\bar{m}_{ii}^{-1}\xi_{i}\intgammah \nablagammah\chi_i\cdot\nablagammah \sum_{j=1}^N\chi_i.
\end{equation}
Since $\Gamma_h$ has no boundary, $\sum_{j=1}^N\chi_i \equiv 1$ and thus
\begin{equation}
\label{expression5heat}
\nablagammah\sum_{j=1}^N\chi_i \equiv 0.
\end{equation}
By combining \eqref{expression1heat}-\eqref{expression5heat}, we finally have
\begin{equation}
-(d\bar{M}^{-1}A\boldxi)_i \leq 0.
\end{equation}
\end{enumerate}
The above points 1. and 2. imply the desired result that $\dot{\xi}_{i}$ is negative. This completes the proof.

\end{proof}
\end{theorem}
\noindent
\begin{theorem}[Maximum principle for \eqref{fullydiscretescheme}]
\label{thm:fdmaximumprinciple}
The fully-discrete solution $\boldxi^n$ with initial data $\boldxi^0$ of scheme \eqref{fullydiscretescheme} fulfills the following maximum principle
\begin{equation}
\label{fullydiscretemaximumprinciple}
0\leq \xi^n_i \leq \max_{\mathbb{R}^N} \left\lbrace \boldxi^0\right\rbrace,\qquad \forall\ i=1,\dots,N,\ \forall\ n\in\mathbb{N},
\end{equation}
if the time step $\tau$ satisfies
\begin{equation}
\label{heattimesteprestriction}
\beta\tau \leq \left(\max_{\boldy\in \Gamma_h}\left\lbrace U^0(\boldy)\right\rbrace\right)^{1-\alpha}.
\end{equation}
In particular, for $\beta=0$, \eqref{fullydiscretemaximumprinciple} holds with no restriction on $\tau$.
\begin{proof}
From the matrix properties \eqref{matrixproperty1} and \eqref{matrixproperty2} we have that, for every $\tau>0$,
\begin{align}
\label{property1again}
&(\bar{M}+d\tau A)^{-1}\bar{M} \geq \vec{0},\\
\label{property2again}
&(\bar{M}+d\tau A)^{-1}\bar{M}\vec{1} = \vec{1}.
\end{align}
In order for the scheme \eqref{fullydiscretescheme} to fulfill the maximum principle \eqref{fullydiscretemaximumprinciple}, it remains to determine a condition on $\tau$ such that
\begin{equation}
\label{positiveness}
\boldxi^n-\tau\beta (\boldxi^n)^\alpha\geq 0,\qquad \forall\ n\in\mathbb{N}.
\end{equation}
Indeed \eqref{fullydiscretescheme}, \eqref{property1again}, \eqref{property2again} and \eqref{positiveness} imply that
\begin{equation}
\label{decreasing}
\max_{\mathbb{R}^N} \left\lbrace\boldxi^{n+1}\right\rbrace \leq \max_{\mathbb{R}^N}\left\lbrace\boldxi^n\right\rbrace, \qquad \forall\ n\in\mathbb{N},
\end{equation}
i.e. the spatial maximum of the fully-discrete solution is not increasing in time. The inequality \eqref{positiveness} may be rewritten elementwise as
\begin{equation}
\label{conditionontau}
\beta\tau \leq (\xi^n_i)^{1-\alpha}, \qquad \forall\ i=1,\dots,N,\ n\in\mathbb{N},
\end{equation}
that is to say
\begin{equation*}
\beta\tau \leq \left(\displaystyle\max_{\mathbb{R}^N} \left\lbrace\boldxi^n\right\rbrace\right)^{1-\alpha} \underset{\eqref{decreasing}}{\leq} \left(\displaystyle\max_{\mathbb{R}^N} \left\lbrace\boldxi^0\right\rbrace\right)^{1-\alpha} = \left(\displaystyle\max_{\boldy\in\Gamma_h} \left\lbrace U^0(\boldy)\right\rbrace\right)^{1-\alpha}
\end{equation*}
which completes the proof.
\end{proof}
\end{theorem}

\section{Reaction-diffusion systems on surfaces}
\label{sec:rdsysyems}
In this section we consider a more general class of surface PDEs that are reaction-diffusion systems of arbitrarily many equations. Analogously to the semilinear parabolic equation \eqref{systemunderstudy}, we apply a lumped finite element space discretization and an IMEX Euler time discretization. We prove that the LSFEM preserves the invariant hyper-rectangles for the semi-discrete and the fully-discrete problems. For the latter case a time step restriction is required.

\subsection{The continuous problem}
\label{sec:continuoussys}
If $\Gamma$ is a compact orientable surface in $\mathbb{R}^3$ without boundary, as in the previous section, and $r\in\mathbb{N}$, let us consider the following reaction-diffusion system of $r$ equations on $\Gamma$:
\begin{equation}
\label{reactiondiffusionsystem}
\begin{cases}
&\dot{u}_1 - d_1\Delta_\Gamma u_1 = f_1(u_1,\dots, u_r),\\
&\hspace{2cm}\vdots\hspace{4cm}(\boldx,t)\in \Gamma\times (0,T],\\
&\dot{u}_r - d_r\Delta_\Gamma u_r = f_r(u_1,\dots, u_r),
\end{cases}
\end{equation}
where $f_1,\dots,f_r$ are $\mathcal{C}^2(\Gamma^r; \mathbb{R})$ reaction kinetics and a $\mathcal{C}^2(\Gamma)$ initial condition is given. As remarked in the Introduction, the following arguments still hold for systems on surfaces with boundary and homogeneous Neumann boundary conditions, i.e. zero \emph{conormal derivative} on $\partial \Gamma$ \cite{acta}. Then, as a special case, planar bounded domains in $\mathbb{R}^2$ with zero-flux boundary conditions could be included in our study. We will confine the present analysis to the case of compact surfaces without boundary to simplify the presentation. In vector form, system \eqref{reactiondiffusionsystem} is given by
\begin{equation}
\label{strongreactiondiffusion}
\begin{cases}
&\dot{\boldu} - D\Delta_\Gamma\boldu = \boldf(\boldu),\qquad (\boldx,t)\in \Gamma\times (0,T],\\
&\boldu(\boldx,0) = \boldu_0(\boldx),\qquad \boldx\in\Gamma,
\end{cases}
\end{equation}
where $D := \text{diag}(d_1,\dots, d_r)$, $\boldu := (u_1,\dots, u_r )^T$, $\Delta_\Gamma\boldu := (\Delta_\Gamma u_1,\dots,\Delta_\Gamma u_r)^T$ and $\boldf (\boldu ) := (f_1(\boldu),\dots, f_r(\boldu))^T$. The weak formulation of \eqref{reactiondiffusionsystem} is: find $u_1,\dots, u_r\in L^2([0,T]; H^1(\Gamma)) \cap L^\infty([0,T]\times \Gamma)$ with $\dot{u}_1,\dots, \dot{u}_r\in L^2([0,T]; H^{-1}(\Gamma))$ such that
\begin{equation}
\label{weakreactiondiffusion}
\begin{cases}
\displaystyle\intgamma \dot{u}_1\varphi_1 + d_1\intgamma \nablagamma u_1\cdot\nablagamma \varphi_1 = \intgamma f_1(\boldu)\varphi_1,\qquad \forall\ \varphi_1\in H^1(\Gamma), \\
\hspace{3cm}\vdots\\
\displaystyle\intgamma \dot{u}_r\varphi_r + d_r\intgamma \nablagamma u_r\cdot\nablagamma \varphi_r = \intgamma f_r(\boldu)\varphi_r,\qquad \forall\ \varphi_r\in H^1(\Gamma).
\end{cases}
\end{equation}
In order to write the corresponding vector formulation we extend all the spatial norms considered throughout the paper to vector-valued functions $\boldv: \Gamma\rightarrow\mathbb{R}^r$ or $\boldV:\Gamma_h\rightarrow\mathbb{R}^r$ as follows. Given a function space $S$, we consider the tensor product norm on $S^r$ defined by
\begin{equation}
\label{tensorproductnorm}
\|\boldv\|_{S^r} := \sqrt{\sum_{i=1}^r \|v_i\|_{S}^2},\qquad \forall \boldv\in S^r.
\end{equation}
For  $p\in [1,+\infty]$, the $L^p([0,T]; S^r)$ norms of space and time dependent functions $\boldu: \Gamma \times [0,T] \rightarrow\mathbb{R}^r$ are defined accordingly. Without any loss of generality, we can write $\|\cdot\|_S$ and $L^p([0,T];S)$ instead of $\|\cdot\|_{S^r}$ and $L^p([0,T];S^r)$, respectively.
Following \cite{barreirathesis}, we introduce the following vector notation:
\begin{equation*}
A:B := \sum_{i=1}^n\sum_{j=1}^m a_{ij}b_{ij},\qquad\forall\ A,B\in \mathbb{R}^{n,m},\quad \forall\ n,m\in\mathbb{N}.
\end{equation*}
We can now write the sum of the equations \eqref{weakreactiondiffusion} as
\begin{equation}
\label{sumweakreactiondiffusion}
\intgamma \dot{\boldu} : \boldvarphi -\intgamma D\nablagamma \boldu : \nablagamma \boldvarphi = \intgamma \boldf(\boldu):\boldvarphi,\qquad \forall\boldvarphi\in (H^1(\Gamma))^r,
\end{equation}
where $\nablagamma \boldu$ is the $r\times 3$ matrix defined by $\nablagamma\boldu := (\nablagamma u_1, \dots, \nablagamma u_r)^T$.

\subsection{Space discretisation}
\label{sec:semi-discretesys}
Analogous to the spatially discretized semilinear parabolic equation \eqref{semi-discreteformulation}, we define the following space discretization for the reaction-diffusion system \eqref{weakreactiondiffusion}: find $U_1,\dots, U_r\in L^2([0,T]; S_h)$ with $\dot{U}_1,\dots, \dot{U}_r\in L^2([0,T]; S_h)$ such that
\begin{equation}
\label{semi-discretereactiondiffusion}
\begin{cases}
\displaystyle\intgamma I_h(\dot{U}_1\phi_1) + d_1\intgamma \nablagamma u_1\cdot\nablagamma \varphi_1 = \intgamma I_h(f_1(\boldU)\phi_1),\qquad \forall \phi_1\in S_h,\\
\hspace{3cm}\vdots\\
\displaystyle\intgamma I_h(\dot{U}_r\phi_r) + d_r\intgamma \nablagamma u_r\cdot\nablagamma \varphi_r = \intgamma I_h(f_r(\boldU)\phi_r),\qquad \forall \phi_r\in S_h.
\end{cases}
\end{equation}
By expressing each component $u_i$ according to \eqref{nodalbasisexpansion}, we have the following matrix form
\begin{equation}
\label{matrixsemi-discretereactiondiffusion}
\begin{cases}
\bar{M}\dot{\boldxi}_1 +d_1A\boldxi_1 = \bar{M}f_1(\boldxi_1,\dots,\boldxi_r),\\
\hspace{2cm}\vdots\\
\bar{M}\dot{\boldxi}_r + d_rA\boldxi_r = \bar{M}f_r(\boldxi_1,\dots, \boldxi_r),
\end{cases}
\end{equation}
where $\bar{M}$ and $A$ are the lumped mass matrix and the stiffness matrix defined in \eqref{lumpedmassmatrix} and \eqref{stiffnessmatrix}, respectively.

\subsection{Time discretization}
\label{sec: fullydiscretesys}
By applying the IMEX Euler method to \eqref{semi-discretereactiondiffusion} we obtain the following fully-discrete method for \eqref{weakreactiondiffusion}: for all $n\in\mathbb{
N}\cup \{0\}$ find $U_1^n,\dots U_r^n \in S_h$ such that
\begin{equation}
\label{fullydiscreteRD}
\begin{cases}
\displaystyle\intgamma I_h\left(\frac{U_1^{n+1}-U_1^n}{\tau}\phi_1^n\right) + d_1\intgamma \nablagamma U_1^{n+1}\cdot\nablagamma \phi_1^n  = \intgamma I_h(f_1(\boldU^n)\phi_1^n),\\
\hspace{3cm}\vdots\\
\displaystyle\intgamma I_h\left(\frac{U_r^{n+1}-U_r^n}{\tau}\phi_r^n\right) + d_r\intgamma \nablagamma U_r^{n+1}\cdot\nablagamma \phi_r^n = \intgamma I_h(f_r(\boldU^n)\phi_r^n), 
\end{cases}
\end{equation}
for all $n\in\mathbb{N}\cup \{0\}$ and $\phi_1^n,\dots, \phi_r^n \in S_h$. We can write the sum of equations \eqref{fullydiscreteRD} as
\begin{equation}
\label{sumfullydiscreteRD}
\intgamma I_h\left(\frac{\boldU^{n+1}-\boldU^n}{\tau}:\boldphi^n\right) + \intgamma D\nablagamma \boldU^{n+1}:\nablagamma \boldphi^n = \intgamma I_h(\boldf(\boldU^n):\boldphi^n),
\end{equation}
for all $n\in\mathbb{N}\cup \{0\}$ and $\boldphi^n \in (S_h)^r$. System \eqref{fullydiscreteRD} can be written in matrix form as
\begin{equation}
\label{fullydiscretereactiondiffusion}
\begin{cases}
\boldxi_1^{n+1} = (\bar{M}+d_1\tau A)^{-1}\bar{M}(\boldxi_1^n + \tau f_1(\boldxi_1^n,\dots,\boldxi_r^n)),\\
\hspace{3cm}\vdots\\
\boldxi_r^{n+1} = (\bar{M}+d_r\tau A)^{-1}\bar{M}(\boldxi_r^n + \tau f_r(\boldxi_1^n,\dots,\boldxi_r^n)),
\end{cases}
\end{equation}
that can be obtained equivalently by applying the IMEX Euler method directly to the ODE system \eqref{matrixsemi-discretereactiondiffusion}.

\subsection{Preservation of the invariant rectangles}
\label{sec:invariantrectangles}
In this section we investigate an interesting property of the lumped finite element discretization of reaction-diffusion systems, which does not hold in the absence of lumping: the preservation of invariant hyper-rectangles. A numerical counterexample will be given in Section \ref{sec:numericaltest}. This preservation property is crucial when the continuous system is known to have an invariant rectangle for two reasons: (i) the solution might be physically meaningless outside a certain range of feasible values, containing the rectangle and (ii) it is a tool to prove stability estimates and error bounds for the semi- and fully-discrete solutions. We recall the following definition given in \cite{smollerbook, taylor1997partial}.
\begin{definition}
\label{def:invariantregion}
For the system \eqref{reactiondiffusionsystem}, a region $\Sigma$ in the phase-space $\mathbb{R}^r$ is said to be a positively invariant region if, whenever the initial condition $\boldu_0$ is in $\Sigma$, $\boldu$ stays in $\Sigma$ as long as it exists and is unique.
\end{definition}

The following theorem has been proven in \cite{smollerbook} when $\Gamma$ is a monodimensional domain in $\mathbb{R}$, in \cite{smollerarticle} when $\Gamma$ is a $k$-dimensional domain in $\mathbb{R}^k$, $k\in\mathbb{N}$ (zero-flux boundary conditions are enforced if the domain is not the whole space) and in \cite{taylor1997partial} in the case in which $\Gamma$ is a Riemannian manifold without boundary. This result provides a sufficient condition for $\Sigma$ to be a positively invariant region in the phase space.
\begin{theorem}[Invariant rectangles for the continuous system \eqref{reactiondiffusionsystem} \cite{taylor1997partial}]
\label{thm:smollertheorem}
Let $\Sigma = \prod_{k=1}^r [m_k,M_k]$ be a hyper-rectangle in the phase space of \eqref{reactiondiffusionsystem}, let $\boldf$ be Lipschitz on $\Sigma$ and let $\mathbf{n}$ be the unit outward vector defined piecewise on $\partial\Sigma$. If
\begin{equation}
\label{smollercondition}
\boldf(\boldu)\cdot\mathbf{n}(\boldu) < 0, \qquad \forall \boldu\in\partial\Sigma,
\end{equation}
then $\Sigma$ is an invariant region for \eqref{reactiondiffusionsystem}.\qed
\end{theorem}
\noindent
Further assumptions on $\boldf$ such that the global existence and uniqueness of the solution are ensured can be found in \cite{taylor1997partial}. We remark that some systems are known to possess an invariant region which do not meet the strict inequality \eqref{smollercondition}. For instance, for many mass-action laws, the positive orthant is invariant \cite{chellaboina2009modeling} even though the flow of $\boldf$ is tangent to this region, instead of strictly inward.

In the following theorems we prove that, under the same assumptions, $\Sigma$ is an invariant region for the semi-discrete \eqref{semi-discretereactiondiffusion} and for the fully-discrete solution \eqref{fullydiscretereactiondiffusion} conditionally on $\tau$, as well. Furthermore, in the fully-discrete case, we will relax the strict inequality \eqref{smollercondition} by requiring non-outward flows, only.

\begin{theorem}[Invariant rectangles for \eqref{semi-discretereactiondiffusion}]
\label{thm:sdinvariantregion}
Let $\Sigma = \prod_{k=1}^r [m_k,M_k]$ be a hyper-rectangle in the phase space, let $\boldf$ be Lipschitz on $\Sigma$ and let $\mathbf{n}$ be the outward unit normal defined piecewise on $\partial \Sigma$. If
\begin{equation}
\label{cond1}
\boldf(\boldu)\cdot\mathbf{n}(\boldu) < 0, \qquad \forall \boldu\in\partial\Sigma,
\end{equation}
then $\Sigma$ is an invariant region for the semi-discrete problem \eqref{semi-discretereactiondiffusion}.
\begin{proof}
We rewrite the semi-discrete problem \eqref{matrixsemi-discretereactiondiffusion} as
\begin{equation}
\label{theoreticalmatrixsemi-discretereactiondiffusion}
\begin{cases}
\dot{\boldxi}_1 +d_1\bar{M}^{-1}A\boldxi_1 = f_1(\boldxi_1,\dots,\boldxi_r),\\
\hspace{2cm}\vdots\\
\dot{\boldxi}_r + d_r\bar{M}^{-1}A\boldxi_r = f_r(\boldxi_1,\dots, \boldxi_r).
\end{cases}
\end{equation}
It suffices to prove that the $rN$-dimensional rectangle $\bar{\Sigma} = \prod_{k=1}^r [m_k,M_k]^N$ is an invariant region for the ODE system \eqref{theoreticalmatrixsemi-discretereactiondiffusion}, i.e. we have to prove that the vector field
\begin{equation*}
\begin{pmatrix}
\dot{\boldxi}_1\\
\vdots\\
\dot{\boldxi}_r
\end{pmatrix} =  -\begin{pmatrix}
d_1\bar{M}^{-1}A & \ &0\\
\ & \ddots & \ \\
0 & \ & d_r\bar{M}^{-1}A
\end{pmatrix} \begin{pmatrix}
{\boldxi}_1\\
\vdots\\
{\boldxi}_r
\end{pmatrix} + \begin{pmatrix}
f_1({\boldxi}_1,\dots,\boldxi_r)\\
\vdots\\
f_r({\boldxi}_1,\dots,\boldxi_r)
\end{pmatrix}
\end{equation*}
computed on every $(rN-1)$-dimensional face of $\bar{\Sigma}$ points toward the interior of $\bar{\Sigma}$. To this end, let $(\boldxi_1,\dots, \boldxi_r)^T$ be a point on $\partial \bar{\Sigma}$. This means that there exist $i=1,\dots,N$ and $k=1,\dots,r$ such that $\xi_{k,i}\in \{m_k,M_k\}$. Suppose $\xi_{k,i} = M_k$; in the case $\xi_{k,i} = m_k$ the proof is analogous. Then
\begin{equation}
\label{boundary}
\xi_{k,j}\leq \xi_{k,i},\quad j\neq i.
\end{equation}
All we have to prove is that $\dot{\xi}_{k,i}$ is negative. To see this, consider that
\begin{enumerate}
\item $f_k(\xi_{1,i},\dots,\xi_{r,i}) = f_k(\xi_{1,i},\dots, M_k,\dots,\xi_{r,i})< 0$ from \eqref{cond1};
\item the $i^{th}$ component of the vector $-d_1\bar{M}^{-1}A\boldxi_k$ is nonpositive. In fact, since $\bar{M}$ is a diagonal matrix, this component is given by
\begin{equation}
\label{expression1}
-(d_k\bar{M}^{-1}A\boldxi_k)_i = -d_k\bar{m}_{ii}^{-1}\sum_{j=1}^N a_{ij}\xi_{k,j}.
\end{equation}
We can split the sum on the right-hand-side by isolating the $a_{ii}\xi_{k,i}$ term:
\begin{equation}
\label{expression2}
-d_k\bar{m}_{ii}^{-1}\sum_{j=1}^N a_{ij}\xi_{k,j} = d_k\bar{m}_{ii}^{-1}\left(-a_{ii}\xi_{k,i} +\sum_{j\in \{1,\dots,N\}\setminus \{i\}} (-a_{ij})\xi_{k,j}\right).
\end{equation}
Since $a_{ij}\leq 0$ for $i\neq j$ from Lemma \ref{delaunaylemma} and $\xi_{k,j}\leq \xi_{k,i}$ for $j\neq i$ from \eqref{boundary}, expression \eqref{expression2} is less than or equal to
\begin{equation}
\label{expression3}
d_k\bar{m}_{ii}^{-1}\xi_{k,i}\left(-a_{ii} +\sum_{j\in\{1,\dots,N\} \setminus \{ i\}} (-a_{ij})\right) = -d_k\bar{m}_{ii}^{-1}\xi_{k,i}\sum_{j=1}^N a_{ij}.
\end{equation}
From the definition of $A$ \eqref{stiffnessmatrix}, expression \eqref{expression3} can be rewritten as
\begin{equation}
\label{expression4}
d_k\bar{m}_{ii}^{-1}\xi_{k,i}\intgammah \nablagammah\chi_i\cdot\nablagammah \sum_{j=1}^N\chi_i.
\end{equation}
Since $\Gamma_h$ has no boundary, $\sum_{j=1}^N\chi_i \equiv 1$ and thus
\begin{equation}
\label{expression5}
\nablagammah\sum_{j=1}^N\chi_i \equiv 0.
\end{equation}
By combining \eqref{expression1}-\eqref{expression5}, we finally have
\begin{equation}
-(d_k\bar{M}^{-1}A\boldxi_k)_i \leq 0.
\end{equation}
\end{enumerate}
These two claims imply the desired fact, i.e. that $\dot{\xi}_{k,i}$ is negative. This completes the proof.
\end{proof}
\end{theorem}

The following theorem is a fully-discrete counterpart of the previous one. Observe that the strictly inward flux condition \eqref{cond1} is replaced by a weaker requirement. This makes the fully-discrete scheme \eqref{fullydiscretereactiondiffusion} somehow more stable than the spatially discrete one \eqref{semi-discretereactiondiffusion}.  The reason for this is that, given a trajectory $u(\boldx,t)$ whose time derivative vanishes at $(\bar{\boldx},\bar{t})$, the function $t\mapsto u(\bar{\boldx},t)$ might still be strictly monotonic, this means that a trajectory may escape a region $\Sigma$ even though the flux of the kinetic is tangent to $\partial\Sigma$.

\begin{theorem}[Invariant rectangles for \eqref{fullydiscretereactiondiffusion}]
\label{thm:fdinvariantregion}
Let $\Sigma = \prod_{k=1}^r [m_k,M_k]$ be a region in the phase space such that
\begin{equation}
\label{cond2}
\boldf(\boldu)\cdot\mathbf{n}(\boldu) \leq 0,\qquad \forall \boldu\in\partial\Sigma.
\end{equation}
For all $k=1,\dots,r$, let $L_k$ be the Lipschitz constant of $f_k$ on $\Sigma$. Then $\Sigma$ is an invariant region for the scheme \eqref{fullydiscretereactiondiffusion} if the time step $\tau$ fulfills
\begin{equation}
\label{fdstabilitycondition}
\tau\leq \frac{1}{\displaystyle\max_{k=1,\dots,r}(L_k)}.
\end{equation}
\begin{proof}
From the matrix properties \eqref{matrixproperty1} and \eqref{matrixproperty2} it follows that, for every $\tau >0$
\begin{align*}
&(\bar{M}+d_k\tau A)^{-1}\bar{M} \geq \vec{0},\qquad \forall\ k=1,\dots,r,\\
&(\bar{M}+d_k\tau A)^{-1}\bar{M}\vec{1} = \vec{1},\qquad \forall k=1,\dots,r.
\end{align*}
In order for the fully-discrete scheme \eqref{fullydiscretereactiondiffusion} to fulfill the theorem, it remains to ensure that
\begin{equation}
\label{cond3}
\begin{split}
m_k \leq \xi^n_{k,i} + \tau f_k(\xi^n_{1,i},\dots,\xi^n_{r,i}) \leq M_k,\\
\forall\ k=1,\dots,r,\ \forall\ i=1,\dots,N,\ \forall\ n\in\mathbb{N}\cup \{0\}.
\end{split}
\end{equation}
Condition \eqref{cond3} is equivalent to
\begin{align}
\label{timebound1}
&\tau \leq \frac{M_k-\xi_{k,i}^n}{f_k(\xi_{1,i}^n,\dots,\xi_{r,i}^n)},\qquad \forall i\ \text{s.t.}\ f_k(\xi_{1,i}^n,\dots,\xi_{r,i}^n)> 0,\\
\label{timebound2}
&\tau \leq \frac{m_k-\xi_{k,i}^n}{f_k(\xi_{i,1}^n,\dots,\xi_{r,i}^n)},\qquad \forall i\ \text{s.t.}\ f_k(\xi_{1,i}^n,\dots,\xi_{r,i}^n)< 0,
\end{align}
for all $k=1,\dots,r$ and $n\in\mathbb{N}$. If $f_k(\xi_{1,i}^n,\dots,\xi_{k,i}^n)>0$, then
\begin{equation}
\label{est1}
f_k(\xi_{1,i}^n,\dots,\xi_{r,i}^n) \leq f_k(M_k) + L_k(M_k-\xi_{k,i}^n) \underset{\eqref{cond2}}{\leq} L_k(M_k-\xi_{k,i}^n).
\end{equation}
If, instead, $f_k(\xi_{1,i}^n,\dots,\xi_{k,i}^n)<0$, then
\begin{equation}
\label{est2}
f_k(\xi_{1,i}^n,\dots,\xi_{r,i}^n) \geq f_k(m_k) - L_k(\xi_{k,i}^n-m_k) \underset{\eqref{cond2}}{\geq} -L_k(\xi_{k,i}^n-m_k).
\end{equation}
Using \eqref{est1} in \eqref{timebound1} and \eqref{est2} in \eqref{timebound2} yields
\begin{equation*}
\tau \leq \frac{1}{L_k},\qquad \forall\ k=1,\dots,r,
\end{equation*}
which completes the proof.
\end{proof}
\end{theorem}
\section{Stability and error analysis}
\label{sec:erroranalysis}
In this section we will prove stability estimates and optimal $L^\infty([0,T],L^2(\Gamma))$ error bounds for the semi-discrete \eqref{semi-discretereactiondiffusion} and the fully-discrete \eqref{fullydiscretereactiondiffusion} solutions of the reaction-diffusion system \eqref{reactiondiffusionsystem} of $r\in\mathbb{N}$ equations. This analysis includes the semilinear parabolic equation \eqref{systemunderstudy}, since \eqref{systemunderstudy} is a special case of the system \eqref{reactiondiffusionsystem} for $r=1$ and $f(u) = -\beta u^\alpha$ and the maximum principle $0\leq u \leq \max u_0$ corresponds to the existence of the invariant region $[0,\max_\Gamma u_0]$. To this end, let us introduce some preliminaries and some basic notations.

The lumped $L^2$ product (see for instance \cite{thomeebook, thomee1980, nochetto1996, garvie2007}) defined by
\begin{equation}
\label{quadraturerule}
(U,V)_h := \int_{\Gamma_h}I_h(UV), \qquad \forall\ U,V\in L^2(\Gamma_h),
\end{equation}
where $I_h$ is given in \eqref{linearinterpolant}, induces the following norm on $S_h$
\begin{equation*}
\|U\|_h = \sqrt{(U,U)_h}, \qquad \forall\ U\in S_h,
\end{equation*}
which is equivalento to $\|\cdot\|_{L^2(\Gamma_h)}$, uniformly with respect to $h$ (see \cite{raviart1973use} for the proof):
\begin{equation}
\label{equivalence}
\|U\|_{L^2(\Gamma_h)} \leq \|U\|_h \leq C\|U\|_{L^2(\Gamma_h)}, \qquad \forall\ U\in S_h,\ \forall h>0.
\end{equation}
\newline
Let us define the "broken" Sobolev space
\begin{equation*}
H^2_h(\Gamma_h) := H^1(\Gamma_h)\cap \displaystyle\prod_{K\in\mathcal{K}_h} H^2(K),
\end{equation*}
endowed with the norm defined by
\begin{equation*}
\|U\|_{H^2_h(\Gamma_h)}^2 := \sum_{K\in\mathcal{K}_h} \|U\|_{H^2(K)}^2,\qquad \forall\ U\in H^2_h(\Gamma_h).
\end{equation*}
For the error in the lumped quadrature rule \eqref{quadraturerule}, if $U\in H^2_h(\Gamma_h)$ and $V\in S_h$, then the following estimate holds (see \cite{thomee1980}):
\begin{equation}
\label{quadratureerror}
|\varepsilon_h (U,V)| := \left|\int_{\Gamma_h} (UV- I_h(UV))\right| \leq ch^2 \|U\|_{H^2_h(\Gamma_h)}\|V\|_{H^1(\Gamma_h)}.
\end{equation} 
Inequalities \eqref{equivalence} and \eqref{quadratureerror} have been proven on planar triangulations in \cite{nochetto1996} and \cite{thomee1980}, respectively. However, since the respective proofs are done piecewise on each triangle, they can be trivially extended to triangulated surfaces with an affine map argument.

The following equivalences between the norms of a function $U$ defined on $\Gamma_h$ and its lifted counterpart $U^\ell$ can be found in \cite{acta}.
\begin{lemma}
Let, $K\in\mathcal{K}_h$, $\tilde{K} := \bolda(T)\subset\Gamma$, where the map $\bolda(\boldx)$ is given in \eqref{onetoone}, and $U:K\rightarrow\mathbb{R}$. If the norms exist, then the following inequalities hold
\begin{align}
\label{equivalence1}
c\|U\|_{L^2(K)} &\leq \|U^\ell\|_{L^2(\tilde{K})} \leq C\|U\|_{L^2(K)};\\
\label{equivalence2}
c\|\nabla_T U\|_{L^2(K)} &\leq \|\nabla_{\tilde{K}} U^\ell\|_{L^2(\tilde{K})} \leq C\|\nabla_K U\|_{L^2(K)};\\
\label{equivalence3}
\|\nabla_K^2 U\|_{L^2(K)} &\leq c(\|\nabla_{\tilde{K}}^2 U^\ell\|_{L^2(\tilde{K})} + h\|\nabla_{\tilde{K}} U^\ell\|_{L^2(\tilde{K})}),
\end{align}
where $\nabla_K^2$ and $\nabla_{\tilde{K}}^2$ are the tangential Hessian on $K$ and $\tilde{K}$, respectively.\qed
\end{lemma}
\noindent
From the previous Lemma we derive the following estimate for the broken $H^2$ norm of $U$.
\begin{lemma}
If $u\in H^2(\Gamma)$, then $u^{-\ell}\in H^2_h(\Gamma_h)$ and
\begin{equation}
\label{equivalence4}
\|u^{-\ell}\|_{H^2_h(\Gamma_h)} \leq C(1+h)\|u\|_{H^2(\Gamma)}.
\end{equation}
\begin{proof}
Let $K\in\mathcal{K}_h$. Then, from \eqref{equivalence1}-\eqref{equivalence3}, we have
\begin{equation}
\label{whydoesthisfollowfromtheabove}
\begin{split}
&\|u^{-\ell}\|_{H^2(K)}^2 = \|u^{-\ell}\|_{L^2(K)}^2 + \|\nabla_K u^{-\ell}\|_{L^2(K)}^2 + \|\nabla_K^2 u^{-\ell}\|_{L^2(K)}^2\\
&\leq \frac{1}{c^2} \|u\|_{L^2(\tilde{K})}^2 + \frac{1}{c^2}\norm{\nabla_{\tilde{K}}u}_{L^2(\tilde{K})}^2\\
&+ c^2\|\nabla^2_{\tilde{K}}u\|_{L^2(\tilde{K})}^2 + c^2h^2\norm{\nabla_{\tilde{K}}u}_{L^2(\tilde{K})}^2\\
&\leq C(1+h^2)\|u\|_{H^2(\tilde{K})}^2.
\end{split}
\end{equation}
Now, from \eqref{whydoesthisfollowfromtheabove}, we have
\begin{equation}
\label{payattention}
\begin{split}
&\|u^{-\ell}\|_{H^2_h(\Gamma_h)}^2 = \sum_{K\in\mathcal{K}_h}\|u^{-\ell}\|_{H^2(K)}^2 \underset{\eqref{whydoesthisfollowfromtheabove}}{\leq} C(1+h^2)\sum_{K\in\mathcal{K}_h}\|u\|_{H^2(\tilde{K})}^2\\
&\leq C(1+h^2)\|u\|_{H^2_h(\Gamma)}^2.
\end{split}
\end{equation}
We remark that, in the last inequality of \eqref{payattention}, the exact equality might not hold, since, being $u^{-\ell}$ only $H^2_h(\Gamma_h)$, its gradient $\nablagammah u^{-\ell}$ might have finite jumps across the edges of the triangulation $\mathcal{K}_h$. This completes the proof.
\end{proof}
\end{lemma}

When lifting integrals, a geometric error must be taken into account. The following equalities hold \emph{(see \cite[p.317]{acta})}
\begin{align}
\label{geometricerror1}
&\intgammah UV \hspace{1.7cm} = \intgamma \frac{U^\ell V^\ell}{\delta_h^\ell},\hspace{2.20cm} \forall\ U,V\in L^2(\Gamma_h),\\
\label{geometricerror2}
&\intgammah\nablagammah U\cdot \nablagammah V = \intgamma \nablagamma U^\ell R_h^T\cdot\nablagamma V^\ell,\qquad\forall\ U,V\in H^1(\Gamma_h),
\end{align}
where $\delta_h^\ell:\Gamma\rightarrow\mathbb{R}$ and $R_h^T:\Gamma\rightarrow\mathbb{R}^{3,3}$ are functions such that (see \cite[p.310]{acta})
\begin{align}
\label{geometricalestimatedelta}
&\left\|1-\frac{1}{\delta_h^\ell}\right\|_{L^\infty(\Gamma)} \leq Ch^2,\\
\label{geometricalestimateR}
&\left\|I-R_h\right\|_{L^\infty(\Gamma)} \leq Ch^2.
\end{align}

For the following proofs we need to define the seminorm $|\cdot|_D$ on $(H^1(\Gamma))^r$ and $(H^1(\Gamma_h))^r$ by
\begin{align}
\label{Dnorm}
&|\boldu|_D^2 := \intgamma D\nablagamma \boldu : \nablagamma \boldu,\qquad \forall\ \boldu\in H^1(\Gamma)^r,\\
\label{Dhnorm}
&|\boldU|_D^2 := \intgammah D\nablagammah \boldU : \nablagammah \boldU,\qquad \forall\ \boldU\in H^1(\Gamma_h)^r,
\end{align}
respectively. Since the diffusion matrix $D$ is diagonal with positive entries, it holds that
\begin{align}
\label{p2}
&\min_{i=1,\dots,r}(d_i) |\boldu|_{H^1(\Gamma)}^2 \leq |\boldu|_D^2 \leq \max_{i=1,\dots,r}(d_i) |\boldu|_{H^1(\Gamma)}^2,\hspace{8mm} \forall \boldu\in (H^1(\Gamma))^r,\\
&\min_{i=1,\dots,r}(d_i) |\boldU|_{H^1(\Gamma_h)}^2 \leq |\boldU|_D^2 \leq \max_{i=1,\dots,r}(d_i) |\boldU|_{H^1(\Gamma_h)}^2,\quad \forall \boldU\in (H^1(\Gamma_h))^r,
\end{align}
i.e. the norms \eqref{Dnorm} and \eqref{Dhnorm} are equivalent to $|\cdot|_{H^1(\Gamma)}$ and $|\cdot|_{H^1(\Gamma_h)}$, respectively.

The following stability estimates are carried out with the usual energy argument. However, thanks to the existence of an invariant region, the estimates will not depend exponentially on time, as the proofs will not rely on Gr\"{o}nwall's lemma. Moreover, we require that the reaction kinetics $\boldf$ in \eqref{strongreactiondiffusion} are Lipschitz only in the invariant region, instead of being globally Lipschitz.
\begin{lemma}[Stability estimates for the weak system \eqref{weakreactiondiffusion}]
If $\boldu$ is the solution of \eqref{weakreactiondiffusion}, $\Sigma = \prod_{k=1}^r [m_k,M_k]$ is an invariant region for \eqref{weakreactiondiffusion}, $\boldf$ is Lipschitz (and thus bounded) on $\Sigma$ and $\boldu_0\in\Sigma$, then
\begin{align}
\label{firststability}
&\sup_{t\in[0,T]} ||\boldu||_{L^2(\Gamma)}^2 + \int_0^T\|\nabla_\Gamma\boldu\|_{L^2(\Gamma)}^2 \leq C\left(T + \|\boldu_0\|_{L^2(\Gamma)}^2\right),\\
\label{secondstability}
&\int_0^T\|\dot{\boldu}\|_{L^2(\Gamma)}^2 + \sup_{t\in [0,T]} \|\nabla_\Gamma \boldu\|_{L^2(\Gamma)}^2 \leq C\left(T + \|\nabla_\Gamma\boldu_0\|_{L^2(\Gamma)}^2\right),
\end{align}
for all $T>0$, where $C$ is a constant independent of $T$ and $\boldu_0$.
\begin{proof}
By setting $\boldvarphi=\boldu$ in \eqref{sumweakreactiondiffusion} we have
\begin{equation}
\label{p1}
\frac{1}{2}\frac{\mathrm{d}}{\mathrm{d}t}\int_\Gamma |\boldu|^2 + |\boldu|_D^2 = \int_\Gamma \boldf(\boldu):\boldu.
\end{equation}
Combining \eqref{p2} and \eqref{p1} we have
\begin{equation*}
\frac{\mathrm{d}}{\mathrm{d}t}\|\boldu\|_{L^2(\Gamma)}^2 + |\boldu|_{H^1(\Gamma)}^2 \leq C \int_\Gamma |\boldf(\boldu):\boldu|.
\end{equation*}
Since $\boldu\in\Sigma$ at all times and $\boldf$ is bounded on $\Sigma$, we obtain
\begin{equation}
\label{ofwhat}
\frac{\mathrm{d}}{\mathrm{d}t}\|\boldu\|_{L^2(\Gamma)}^2 + |\boldu|_{H^1(\Gamma)}^2 \leq C.
\end{equation}
By integrating both sides of \eqref{ofwhat} over $[0,T]$, estimate \eqref{firststability} follows.

To prove the second estimate, we set $\boldvarphi = \dot{\boldu}$ in \eqref{sumweakreactiondiffusion} and obtain:
\begin{equation}
\label{p3}
\int_\Gamma |\dot{\boldu}|^2 + \frac{1}{2}\frac{\mathrm{d}}{\mathrm{d}t}\int_\Gamma D\nabla_\Gamma\boldu:\nabla_\Gamma\boldu \leq \int_\Gamma |\boldf(\boldu)||\dot{\boldu}|,
\end{equation}
but, since $\boldf$ is bounded on $\Sigma$, we have
\begin{equation}
\label{p4}
\int_\Gamma |\boldf(\boldu)||\dot{\boldu}| \leq \frac{1}{2}\int_\Gamma |\boldf(\boldu)|^2 + \frac{1}{2}\int_\Gamma |\dot{\boldu}|^2 \leq C + \frac{1}{2}\int_\Gamma |\dot{\boldu}|.
\end{equation}
Combining \eqref{p3} and \eqref{p4} we have
\begin{equation*}
\|\dot{\boldu}\|_{L^2(\Gamma)}^2 + \frac{\mathrm{d}}{\mathrm{d}t}|\boldu|_D^2 \leq C,
\end{equation*}
from which, by integrating on $[0,T]$ we obtain
\begin{equation}
\label{p5}
\int_0^T\|\dot{\boldu}\|_{L^2(\Gamma)}^2 + |\boldu|_D^2 \leq CT + |\boldu_0|_D^2.
\end{equation}
Combining \eqref{p5} with \eqref{p2}, we have
\begin{equation*}
\int_0^T\|\dot{\boldu}\|_{L^2(\Gamma)}^2 + |\boldu|_{H^1(\Gamma)}^2 \leq C\left(T + |\boldu_0|_{H^1(\Gamma)}^2\right),
\end{equation*}
from which we obtain estimate \eqref{secondstability}. \qed
\end{proof}
\end{lemma}
\noindent
The following lemmas show analogous estimates for the semi- and fully-discrete problems.

\begin{lemma}[Stability estimates for the semi-discrete system \eqref{semi-discretereactiondiffusion}]
If $\boldU$ is the solution of \eqref{semi-discretereactiondiffusion}, $\Sigma = \prod_{k=1}^r [m_r,M_r]$ is an invariant region for \eqref{semi-discretereactiondiffusion}, $\boldf$ is Lipschitz on $\Sigma$ and $\boldU_0\in\Sigma$, then
\begin{align}
\label{firststabilityestimate}
&\sup_{t\in[0,T]} ||\boldU||_{L^2(\Gamma_h)}^2 + \int_0^T\|\nabla_\Gamma\boldU\|_{L^2(\Gamma_h)}^2 \leq C\left(T + \|\boldU_0\|_{L^2(\Gamma_h)}^2\right),\\
\label{secondstabilityestimate}
&\int_0^T\|\dot{\boldU}\|_{L^2(\Gamma_h)}^2 + \sup_{t\in [0,T]} \|\nabla_\Gamma \boldU\|_{L^2(\Gamma_h)}^2 \leq C\left(T + \|\nabla_\Gamma\boldU_0\|_{L^2(\Gamma_h)}^2\right),
\end{align}
for all $T>0$, where $C$ is a constant independent of $T$ and $\boldU_0$.
\begin{proof}
We proceed exactly as in the previous Lemma in order to obtain analogous estimates in the norm $\|\cdot\|_h$ and then we use the equivalence \eqref{equivalence} between the norms $\|\cdot\|_h$ and $\|\cdot\|_{L^2(\Gamma_h)}$ on $S_h$, uniformly in $h$.
\end{proof}
\end{lemma}

\begin{lemma}[Stability estimates for the fully-discrete system \eqref{fullydiscreteRD}]
Let $\tau>0$. If $\boldU^i$, $i=0,\dots, \frac{T}{\tau}$, is the solution of \eqref{fullydiscretereactiondiffusion}, $\Sigma = \prod_{k=1}^r [m_r,M_r]$ is an invariant region for \eqref{fullydiscretereactiondiffusion}, $\boldf$ is Lipschitz on $\Sigma$ and $\boldU_0\in\Sigma$, then
\begin{align}
\label{firstfullydiscretestability}
&\|\boldU^{n+1}\|_{L^2(\Gamma_h)}^2 + \tau\sum_{i=0}^n \|\nablagammah \boldU^{i+1}\|_{L^2(\Gamma_h)}^2 \leq C(\|\boldU^0\|_{L^2(\Gamma_h)}+T),\\
\label{secondfullydiscretestability}
\frac{1}{\tau}\sum_{i=0}^n &\|\boldU^{i+1}-\boldU^i\|_{L^2(\Gamma_h)}^2 + \|\nabla_{\Gamma_h} \boldU^{n+1}\|_{L^2(\Gamma_h)}^2 \leq C(\|\nablagammah \boldU^0\|_{L^2(\Gamma_h)}^2 + T),
\end{align}
for all $n =1,\dots, \frac{T}{\tau}$ and $T>0$, where $C$ is a constant independent of $T$ and $\boldU_0$.
\begin{proof}
By testing \eqref{sumfullydiscreteRD} with $\boldphi^i = \boldU^{i+1}$ we have
\begin{equation*}
\frac{1}{\tau}\left(\|\boldU^{i+1}\|_h^2 - \intgammah I_h(\boldU^{i}:\boldU^{i+1})\right) + |\boldU^{i+1}|_D^2 = \intgammah I_h(\boldf(\boldU^i):\boldU^{i+1}).
\end{equation*}
After multiplying by $\tau$, Cauchy-Schwarz inequality yields
\begin{equation*}
\|\boldU^{i+1}\|_h^2 + \tau|\boldU^{i+1}|_D \leq \|\boldU^{i+1}\|_h\|\boldU^i\|_h + \tau\|\boldf(\boldU)^i\|_h\|\boldU^{i+1}\|_h.
\end{equation*}
Since $\boldU^i$ and $\boldU^{i+1}\in\Sigma$ and $\boldf$ is Lipschitz on $\Sigma$, the last term on the right hand side is bounded by some constant $C>0$:
\begin{equation*}
\|\boldU^{i+1}\|_h^2 + \tau|\boldU^{i+1}|_D \leq \|\boldU^{i+1}\|_h\|\boldU^i\|_h + C\tau.
\end{equation*}
Young's inequality yields
\begin{equation*}
\|\boldU^{i+1}\|_h^2 + \tau|\boldU^{i+1}|_D^2 \leq \|\boldU^i\|_h^2 + C\tau.
\end{equation*}
We sum for $i=0,\dots,n$ to obtain
\begin{equation*}
\|\boldU^{n+1}\|_h^2 + \tau\sum_{i=0}^n |\boldU^{i+1}|_D^2 \leq \|\boldU^0\|_h^2 + Cn\tau.
\end{equation*}
By using \eqref{equivalence}, the equivalence between $|\cdot|_D$ and $|\cdot|_{H^1(\Gamma_h)}$ and $n=1,\dots,\frac{T}{\tau}$, \eqref{firstfullydiscretestability} follows immediately.

By testing \eqref{sumfullydiscreteRD} with $\boldphi^i=\boldU^{i+1}-\boldU^i$ we have
\begin{equation*}
\begin{split}
\frac{1}{\tau}\|\boldU^{i+1}-\boldU^i\|_h^2 + |\boldU^{i+1}|_D^2 & - \intgammah D\nablagammah \boldU^{i+1}:\nablagammah \boldU^{i}\\
&= \intgammah I_h(\boldf(\boldU^i):(\boldU^{i+1}-\boldU^i)).
\end{split}
\end{equation*}
Cauchy-Schwarz inequality yields
\begin{equation*}
\frac{1}{\tau}\|\boldU^{i+1}-\boldU^i\|_h^2 + |\boldU^{i+1}|_D^2 \leq |\boldU^{i+1}|_D|\boldU^i|_D+ \|\boldf(\boldU^{i})\|_h\|\boldU^{i+1}-\boldU^i\|_h.
\end{equation*}
Since $\boldf$ is Lipschitz -and thus bounded- on $\Sigma$, say $\max_\Sigma \boldf = C$, we can bound the last term in the right hand side as follows:
\begin{equation*}
\begin{split}
\frac{1}{\tau}\|\boldU^{i+1}-\boldU^i\|_h^2 + |\boldU^{i+1}|_D^2 \leq &|\boldU^{i+1}|_D|\boldU^i|_D + C\|\boldU^{i}\|_h\|\boldU^{i+1}-\boldU^i\|_h.
\end{split}
\end{equation*}
Young's inequality yields
\begin{equation*}
\begin{split}
\frac{1}{\tau}\|\boldU^{i+1}-\boldU^i\|_h^2 + |\boldU^{i+1}|_D^2 \leq & \frac{1}{2}(| \boldU^i|_D^2 + | \boldU^{i+1}|_D^2) + C\tau\\
 + &\frac{1}{2\tau}\|\boldU^{i+1}-\boldU^i\|_h^2.
\end{split}
\end{equation*}
Rearranging terms and multiplying by $2$ we have
\begin{equation}
\label{tobesummed}
\frac{1}{\tau}\|\boldU^{i+1}-\boldU^i\|_h^2 + | \boldU^{i+1}|_D^2 \leq |\boldU^i|_D^2 + C\tau.
\end{equation}
By summing \eqref{tobesummed} for $i=0,\dots,n$ we have
\begin{equation*}
\frac{1}{\tau}\sum_{i=0}^n\| \boldU^{i+1}-\boldU^i\|_h^2 + |\boldU^{n+1}|_D^2 \leq |\boldU^0|_D^2 + Cn\tau.
\end{equation*}
By using \eqref{equivalence}, the equivalence between $|\cdot|_D$ and $|\cdot|_{H^1(\Gamma_h)}$ and $n=1,\dots,\frac{T}{\tau}$, \eqref{secondfullydiscretestability} finally follows.
\end{proof}
\end{lemma}

To prove the convergence of the semi- and fully-discrete methods, we will adopt the surface Ritz projection considered in \cite{du2011finite, elliott2015evolving, lubich2015variational}.
\begin{definition}
Given $u:[0,T]\rightarrow H^1(\Gamma)$, the Ritz projection of $u$ is the unique function $\bar{U}: [0,T]\rightarrow S_h$ such that
\begin{equation}
\label{ritzdefinition}
\begin{split}
&\intgammah \nablagammah \bar{U} \cdot \nablagammah \varphi = \intgamma \nablagamma u\cdot \nablagamma \varphi^\ell,\qquad \forall\ \varphi\in S_h,\\
&\intgammah \bar{U} = \intgamma u.
\end{split}
\end{equation}
\end{definition}
\noindent
We remark that this definition is different from the one considered in \cite{elliott2013}. The following error estimates for the Ritz projection can be found in \cite{du2011finite, elliott2015evolving}.
\begin{theorem}[Error estimates for the Ritz projection]
\label{ritztheorem}
Given $u:[0,T] \rightarrow H^2(\Gamma)$ such that $\dot{u} : [0,T]\rightarrow H^2(\Gamma)$, the error in the Ritz projection satisfies the following bounds
\begin{align}
\label{firstritzestimate}
&\|u-\bar{U}^\ell\|_{L^2(\Gamma)} + h \|\nabla_\Gamma(u-\bar{U}^\ell)\|_{L^2(\Gamma)} \leq ch^2\|u\|_{H^2(\Gamma)},\\
\label{secondritzestimate}
&\|\dot{u}-\dot{\bar{U}}^\ell\|_{L^2(\Gamma)} + h \|\nabla_\Gamma (\dot{u}-\dot{\bar{U}}^\ell)\|_{L^2(\Gamma)} \leq ch^2(\|u\|_{H^2(\Gamma)} + \|\dot{u}\|_{H^2(\Gamma)}).
\end{align} \qed
\end{theorem}
\noindent
If $\boldu$ is a vector function, we will denote with $\bar{\boldU}$ its componentwise Ritz projection and the estimates \eqref{firstritzestimate}-\eqref{secondritzestimate} still hold in the tensor product norms \eqref{tensorproductnorm}. An $L^\infty([0,T],L^2(\Gamma))$ error bound for the semi-discrete solution has been proven in \cite{thomee1980} on planar domains. Here we extend this result to triangulated surfaces.

\begin{theorem}[Error estimate for the semi-discrete solution \eqref{semi-discretereactiondiffusion}]
\label{thm:SDerroranalysis}
Assume that $\Sigma$ is an invariant region for \eqref{weakreactiondiffusion} and \eqref{semi-discretereactiondiffusion}, that $\boldf\in\mathcal{C}^2(\Sigma)$ and that $\boldu_0,\boldU_0\in\Sigma$. If the solution $\boldu$ of \eqref{weakreactiondiffusion} and its time derivative $\dot{\boldu}$ are $L^{\infty}([0,T]; H^2(\Gamma))$ and $\|\boldu_0-\boldU_0^\ell\|_{L^2(\Gamma)} \leq ch^2$, then the following estimate holds
\begin{equation}
\|\boldu- \boldU^\ell\|_{L^2(\Gamma)} \leq C(\boldu,T) h^2,
\end{equation}
where $C(\boldu,T)$ is a constant depending on $\boldu$ and $T$.
\begin{proof}
Let us write the error as
\begin{equation}
\label{spliterrorSD}
\boldU^{\ell}-\boldu = (\boldU^{\ell}-{\bar{\boldU}}^\ell) + ({\bar{\boldU}}^\ell- {\boldu}) =: {\boldtheta}^\ell + {\boldrho}^\ell.
\end{equation}
Since $\boldu$ and $\dot{\boldu}$ are $L^{\infty}([0,T],H^2(\Gamma))$, from the error estimates \eqref{firstritzestimate}-\eqref{secondritzestimate} for the Ritz projection and \eqref{equivalence1}-\eqref{equivalence2} we have that
\begin{align}
\label{ensuresthatSD}
&\norm{\boldrho}_{L^2(\Gamma_h)} \leq C\|\boldrho^\ell\|_{L^2(\Gamma)} = C\|\bar{\boldU}^\ell-\boldu\|_{L^2(\Gamma)} \leq Ch^2\norm{\boldu}_{H^2(\Gamma)},\\
\label{dotensuresthatSD}
&\norm{\dot{\boldrho}}_{L^2(\Gamma_h)} + h \|\nablagammah\dot{\boldrho}\|_{L^2(\Gamma_h)}\leq Ch^2(\norm{\boldu}_{H^2(\Gamma)} + \norm{\dot{\boldu}}_{H^2(\Gamma)}).
\end{align}

It remains to show the convergence for $\boldtheta^\ell$ in \eqref{spliterrorSD}. For the sake of simplicity, we derive an estimate for $\boldtheta$ in the norm $\|\cdot\|_h$ and then we will use \eqref{equivalence} and \eqref{equivalence1} to estimate $\|\boldtheta^\ell\|_{L^2(\Gamma)}$. The continuous problem \eqref{weakreactiondiffusion}, the semi-discrete formulation \eqref{semi-discretereactiondiffusion}, the definition of Ritz projection \eqref{ritzdefinition} and the relations \eqref{geometricerror1}, \eqref{geometricerror2}, imply that
\begin{equation}
\label{p0SD}
\begin{split}
&\intgammah I_h(\dot{\boldtheta}:\boldphi) + \intgammah D\nablagammah \boldtheta: \nablagammah\boldphi =\intgammah I_h((\boldf(\boldU)-\boldf(\boldu^{-\ell})):\boldphi)\\
& + \varepsilon_h(\boldf(\boldu^{-\ell}),\boldphi) +\intgamma \left( 1-\frac{1}{\delta_h^\ell}\right) \boldf(\boldu):\boldphi^{\ell} -\intgammah \dot{\boldrho}:\boldphi + \varepsilon_h(\dot{\boldrho},\boldphi)\\
& +\intgamma \left( 1-\frac{1}{\delta_h^\ell}\right) \boldu:\boldphi^{\ell}.
\end{split}
\end{equation}
In \eqref{p0SD} we choose $\boldphi=\boldtheta$. For the first term of \eqref{p0SD} we observe that
\begin{equation}
\label{p0SDleft}
\intgammah I_h(\dot{\boldtheta} : \boldtheta) =\frac{1}{2}\frac{\mathrm{d}}{\mathrm{d}t}\norm{\boldtheta}_h^2.
\end{equation}
We estimate the single terms on the right hand side of \eqref{p0SD} in turn. By using the Cauchy-Schwarz inequality, the Lipschitz continuity of $\boldf$, the definition of $\boldtheta$, \eqref{equivalence}, \eqref{equivalence1} and \eqref{ensuresthatSD}, we have that
\begin{equation}
\label{p1SD}
\begin{split}
&\semin{\intgammah I_h((\boldf(\boldU)-\boldf(\boldu^{-\ell})) : \boldtheta)} \leq \|\boldf(\boldU)-\boldf(\boldu^{-\ell})\|_h \norm{\boldtheta}_h\\
& \leq C\|\boldU-\boldu^{-\ell}\|_h\norm{\boldtheta}_h \leq C\left(\|\boldrho\|_{L^2(\Gamma)}+ \norm{\boldtheta}_h\right)\norm{\boldtheta}_h\\
&= C(\boldu)(h^2+\norm{\boldtheta}_h) \norm{\boldtheta}_h.
\end{split}
\end{equation}
By using the estimate \eqref{quadratureerror} for $\varepsilon_h$, \eqref{equivalence4}, the regularity assumptions $\boldf\in \mathcal{C}^2(\Sigma)$ and $\boldu\in L^{\infty}([0,T], H^2(\Gamma))$, and by applying the chain rule to the composite function $\boldf(\boldu)$ it follows that
\begin{equation}
\label{p2SD}
\begin{split}
&\semin{\varepsilon_h(\boldf(\boldu^{-\ell}),\boldtheta) } \leq Ch^2\|\boldf(\boldu^{-\ell})\|_{H^2_h(\Gamma_h)}\norm{\boldtheta}_{H^1(\Gamma_h)} \leq\\
& C(1+h)h^2\|\boldf(\boldu)\|_{H^2(\Gamma)}\norm{\boldtheta}_{H^1(\Gamma_h)} \leq \\
&C(h^2+h^3)\|\boldf\|_{\mathcal{C}^2(\Sigma)}\|\boldu\|_{H^2(\Gamma)}\norm{\boldtheta}_{H^1(\Gamma)} \leq C(h^2+h^3)\norm{\boldtheta}_{H^1(\Gamma_h)}.
\end{split}
\end{equation}
Since $\boldf$ is Lipschitz over the compact region $\Sigma$, then $\boldf\in L^\infty(\Sigma)$. Hence, by using the Cauchy-Schwarz inequality and the geometric estimate \eqref{geometricalestimatedelta} we have
\begin{equation}
\label{p3SD}
\begin{split}
&\left|\intgamma \left( 1-\frac{1}{\delta_h^\ell}\right) \boldf(\boldu):\boldtheta\right| \leq \left\| 1-\frac{1}{\delta_h^\ell}\right\|_{L^\infty(\Gamma)}\|\boldf(\boldu)\|_{L^2(\Gamma)}\|\boldtheta\|_{L^2(\Gamma)}\\
& \leq Ch^2\|\boldtheta\|_{L^2(\Gamma)}.
\end{split}
\end{equation}
From the Cauchy-Schwarz inequality, the error estimate \eqref{dotensuresthatSD} and \eqref{equivalence1} we have
\begin{equation}
\label{p6SD}
\quad \left|\intgammah \dot{\boldrho}:\boldtheta\right| \leq C\|\dot{\boldrho}\|_{L^2(\Gamma_h)}\|\boldtheta\|_{L^2(\Gamma_h)} \leq C(\boldu)h^2\|\boldtheta\|_{L^2(\Gamma_h)}.
\end{equation}
From the estimate \eqref{quadratureerror} for $\varepsilon_h$, the estimate \eqref{dotensuresthatSD} for $\boldrho$, \eqref{equivalence1},\eqref{equivalence2} and \eqref{equivalence4} we have
\begin{equation}
\label{p7SD}
\begin{split}
&\semin{\varepsilon_h(\dot{\boldrho},\boldtheta)} \leq Ch^2\|\dot{\boldrho}\|_{H^2_h(\Gamma_h)}\norm{\boldtheta}_{H^1(\Gamma_h)}\\
&= Ch^2(\|\dot{\boldrho}\|_{H^1(\Gamma_h)} + |\dot{\boldrho}|_{H^2_h(\Gamma_h)}) \norm{\boldtheta}_{H^1(\Gamma_h)}\\
&=  Ch^2(\|\dot{\boldrho}\|_{H^1(\Gamma_h)} + |\dot{\boldu}^{-\ell}|_{H^2_h(\Gamma)})\norm{\boldtheta}_{H^1(\Gamma_h)} \leq \\
& Ch^2(C(\boldu)h + (1+h)\|\boldu\|_{H^2(\Gamma)})\norm{\boldtheta}_{H^1(\Gamma_h)}\leq C(\boldu)(h^2+h^3) \norm{\boldtheta}_{H^1(\Gamma_h)},
\end{split}
\end{equation}
where $|\cdot|_{H^2_h(\Gamma_h)}$ denotes the broken $H^2$ seminorm on $\Gamma_h$. Cauchy-Schwarz inequality, \eqref{equivalence1}, the geometric estimate \eqref{geometricalestimatedelta} and the stability bound \eqref{firststability} yield
\begin{equation}
\label{p9SD}
\begin{split}
&\left|\intgamma \left( 1-\frac{1}{\delta_h^\ell}\right) \boldu:\boldtheta^\ell\right| \leq \left\|1-\frac{1}{\delta_h^\ell}\right\|_{L^\infty(\Gamma)}\|\boldu\|_{L^2(\Gamma)}\|\boldtheta\|_{L^2(\Gamma_h)}\\
& \leq Ch^2\|\boldtheta\|_{L^2(\Gamma_h)}.
\end{split}
\end{equation}
Combining \eqref{p0SD}--\eqref{p9SD}, using \eqref{equivalence}, \eqref{equivalence1} and \eqref{equivalence2}, we have
\begin{equation}
\label{p10SD}
\begin{split}
&\frac{1}{2}\frac{\mathrm{d}}{\mathrm{d}t}\norm{\boldtheta}_h^2 + m\norm{\nablagammah \boldtheta}_{L^2(\Gamma_h)} \leq C(\boldu)(h^2+h^3+\norm{\boldtheta}_h)\norm{\boldtheta}_{H^1(\Gamma_h)}\\
&\leq C(\boldu,m)(h^4+h^6+\norm{\boldtheta}_h^2) + m\norm{\boldtheta}_{H^1(\Gamma_h)}^2,
\end{split}
\end{equation}
where $m=\min_{k=1,\dots,r}\{d_k\}$. Canceling $m\|\nablagammah\boldtheta\|_{L^2(\Gamma_h)}$ on both sides of \eqref{p10SD}, and again using \eqref{equivalence}, we have that
\begin{equation*}
\frac{\mathrm{d}}{\mathrm{d}t}\norm{\boldtheta}_h^2 \leq C(\boldu)(h^4+h^6) + C(\boldu)\norm{\boldtheta}_h^2.
\end{equation*}
Using Gr\"{o}nwall's lemma, the assumption $\|\boldtheta_0^\ell\|_{L^2(\Gamma)}\leq Ch^2$, \eqref{equivalence} and \eqref{equivalence1}, we obtain
\begin{equation*}
\|\boldtheta^\ell\|_{L^2(\Gamma)}^2 \leq C(\boldu,T)(h^4+h^6),
\end{equation*}
which yields the desired result.
\end{proof}
\end{theorem}

\noindent In a similar fashion, following the approach in \cite{lakkis2013implicit} and \cite{thomee1980}, we obtain the following $L^\infty([0,T],L^2(\Gamma))$ error estimate for the fully-discrete solution \eqref{fullydiscreteRD}.

\begin{theorem}[Error estimate for the fully-discrete solution \eqref{fullydiscreteRD}]
\label{thm: FDerroranalysis}
Assume that $\Sigma$ is an invariant region for \eqref{weakreactiondiffusion} and \eqref{fullydiscreteRD}, that $\boldf\in\mathcal{C}^2(\Sigma)$ and that $\boldu_0,\boldU_0\in\Sigma$. If the solution $\boldu$ of \eqref{weakreactiondiffusion} and its time derivative $\dot{\boldu}$ are $L^{\infty}([0,T]; H^2(\Gamma))$, $\ddot{\boldu}$ is $L^\infty([0,T]; L^2(\Gamma))$ and $\|\boldu_0-\boldU_0^\ell\|_{L^2(\Gamma)} \leq ch^2$, then the following estimate holds
\begin{equation}
\|\boldu^n - \boldU^{\ell,n}\|_{L^2(\Gamma)} \leq C(\boldu,T)(h^2+h^3+\tau),
\end{equation}
where $\boldu^n$ is the exact solution at time $t_n:=n\tau$ and $C(\boldu,T)$ is a constant depending on $\boldu$ and $T$.
\begin{proof}
Let us write the error as
\begin{equation}
\label{spliterrorFD}
\boldU^{\ell,n}-\boldu^n = (\boldU^{\ell,n}-{\bar{\boldU}}^{\ell,n}) + ({\bar{\boldU}}^{\ell,n} - {\boldu}^n) =: {\boldtheta}^{\ell,n} + {\boldrho}^{\ell,n},\quad  \forall n,
\end{equation}
and the discrete time derivative of any function $\boldphi: \Gamma_h\times [0,T]\rightarrow\mathbb{R}^r$ as
\begin{equation*}
\bar{\partial}\boldphi^n := \frac{\boldphi^n-\boldphi^{n-1}}{\tau},\quad \forall n. 
\end{equation*}
Since $\boldu$ and $\dot{\boldu}$ are $L^{\infty}([0,T],H^2(\Gamma))$, from \eqref{equivalence1}, \eqref{equivalence2}, \eqref{firstritzestimate} and \eqref{secondritzestimate}, we have that
\begin{align}
\label{ensuresthatFD}
&\|\boldrho^n\|_{L^2(\Gamma_h)} \leq C\|\boldrho^{\ell,n}\|_{L^2(\Gamma)} = \|\bar{\boldU}^{\ell,n}-\boldu^n\|_{L^2(\Gamma)} \leq ch^2\norm{\boldu^n}_{H^2(\Gamma)},\quad \forall n,\\
\label{dotensuresthatFD}
&\|\dot{\boldrho}^n\|_{L^2(\Gamma_h)} + h \|\nablagammah\dot{\boldrho}^n\|_{L^2(\Gamma_h)} \leq ch^2(\|\boldu^n\|_{H^2(\Gamma)} + \|\dot{\boldu}^n\|_{H^2(\Gamma)}),\quad \forall n.
\end{align}
It remains to show the convergence for $\boldtheta^{\ell,n}$ in \eqref{spliterrorFD}. To this end, we derive an estimate for $\boldtheta^{n}$ in the $L^2(\Gamma_h)$ norm and then use \eqref{equivalence} and \eqref{equivalence1} to estimate $\|\boldtheta^{\ell,n}\|_{L^2(\Gamma)}$.

The continuous problem \eqref{weakreactiondiffusion} and the fully-discrete formulation \eqref{fullydiscreteRD}, the definition of Ritz projection \eqref{ritzdefinition}, the relations \eqref{geometricerror1} and \eqref{geometricerror2}, imply that
\begin{equation}
\label{p1fully}
\begin{split}
&\intgammah I_h(\bar{\partial}\boldtheta^{n}:\boldphi^n) + \intgammah D\nablagammah \boldtheta^{n}: \nablagammah\boldphi^n = \varepsilon_h(\boldf(\boldu^{-\ell,n-1}),\boldphi^n)\\
+&\intgammah I_h((\boldf(\boldU^{n-1})-\boldf(\boldu^{-\ell,n-1})):\boldphi^n)+\intgamma \left( 1-\frac{1}{\delta_h^\ell}\right) \boldf(\boldu^{n-1}):\boldphi^{\ell,n}\\
+&\intgamma (\boldf(\boldu^{n-1})-\boldf(\boldu^n)):\boldphi^{\ell,n}-\intgammah \bar{\partial}\boldrho^{n}:\boldphi^n + \varepsilon_h(\bar{\partial}\boldrho^{n},\boldphi^n)\\
-&\intgammah (\bar{\partial}-\partial_t)\boldu^{-\ell,n}:\boldphi^n +\intgamma \left( 1-\frac{1}{\delta_h^\ell}\right) \boldu^n:\boldphi^{\ell,n}.
\end{split}
\end{equation}
In \eqref{p1fully} we choose $\boldphi^n=\boldtheta^{n}$. For the first term in \eqref{p1fully} we observe that, from Young's inequality,
\begin{equation}
\label{p1fullyleft}
\intgammah I_h\left(\bar{\partial}\boldtheta^{n} : \boldtheta^{n}\right) \geq \frac{1}{2\tau}( \norm{\boldtheta^{n}}_h^2 - \|\boldtheta^{n-1}\|_h^2).
\end{equation}
We estimate the single terms on the right hand side of \eqref{p1fully} in turn. From the Cauchy-Schwarz inequality, the Lipschitz continuity of $\boldf$, the definition of $\boldtheta^{n}$, \eqref{equivalence} and \eqref{ensuresthatFD}, it follows that
\begin{equation}
\label{p1FD}
\begin{split}
&\semin{\intgammah I_h((\boldf(\boldU^{n-1})-\boldf(\boldu^{-\ell,n-1})) : \boldtheta^{n})}\\
& \leq \|\boldf(\boldU^{n-1})-\boldf(\boldu^{-\ell,n-1})\|_h \|\boldtheta^{n}\|_h\leq C\|\boldU^{n-1}-\boldu^{-\ell,n-1}\|_h\|\boldtheta^{n}\|_h\\
& \leq C(\|\boldrho^{n-1}\|_{L^2(\Gamma)}+ \|\boldtheta^{n-1}\|_h)\|\boldtheta^{n}\|_h\leq C(\boldu )(h^2 +\|\boldtheta^{n-1}\|_h) \|\boldtheta^{n}\|_h.
\end{split}
\end{equation}
From the estimate \eqref{quadratureerror} for $\varepsilon_h$ and \eqref{equivalence4}, it follows that
\begin{equation}
\label{p2FD}
\begin{split}
&\semin{\varepsilon_h(\boldf(\boldu^{-\ell,n-1}),\boldtheta^{n})} \leq Ch^2\|\boldf(\boldu^{-\ell,n-1})\|_{H^2_h(\Gamma_h)}\|\boldtheta^{n}\|_{H^1(\Gamma_h)}\\
&\leq C(1+h)h^2\|\boldf(\boldu^{n-1})\|_{H^2(\Gamma)}\|\boldtheta^{n}\|_{H^1(\Gamma_h)}\\
& \leq C(1+h)h^2\|\boldf\|_{\mathcal{C}^2(\Sigma)}\|\boldu^{n-1}\|_{H^2(\Gamma)}\|\boldtheta^{n}\|_{H^1(\Gamma_h)}\\
&\leq C(h^2+h^3)\|\boldtheta^{n}\|_{H^1(\Gamma_h)},
\end{split}
\end{equation}
where we have exploited the regularity assumptions $\boldf\in \mathcal{C}^2(\Sigma)$ and $\boldu\in L^{\infty}([0,T], H^2(\Gamma))$. Since $\boldf$ is Lipschitz over the compact region $\Sigma$ then $\boldf\in L^\infty(\Sigma)$. This fact, together with Cauchy-Schwarz inequality, \eqref{equivalence1} and the geometric estimate \eqref{geometricalestimatedelta}, yields
\begin{equation}
\label{p3FD}
\begin{split}
&\semin{\intgamma \left( 1-\frac{1}{\delta_h^\ell}\right) \boldf(\boldu^{n-1}):\boldtheta^{\ell,n}}\\
& \leq \left\| 1-\frac{1}{\delta_h^\ell}\right\|_{L^\infty(\Gamma)}\|\boldf(\boldu^{n-1})\|_{L^2(\Gamma_h)}\|\boldtheta^n\|_{L^2(\Gamma)} \leq Ch^2\|\boldtheta^n\|_{L^2(\Gamma_h)}.
\end{split}
\end{equation}
Cauchy-Schwarz inequality yields, together with \eqref{equivalence1} and the stability estimate \eqref{secondstability},
\begin{equation}
\label{p4FD}
\begin{split}
&\semin{\intgamma (\boldf(\boldu^{n-1})-\boldf(\boldu^n)):\boldtheta^{\ell,n}} \leq \|\boldf(\boldu^{n-1})-\boldf(\boldu^n)\|_{L^2(\Gamma)}\norm{\boldtheta^n}_{L^2(\Gamma_h)}\\
& \leq C\|\boldu^n-\boldu^{n-1}\|_{L^2(\Gamma)}\norm{\boldtheta^n}_{L^2(\Gamma_h)} = \norm{ \int_{t_{n-1}}^{t_n}\dot{\boldu}}_{L^2(\Gamma)}\|\boldtheta^n\|_{L^2(\Gamma_h)}\\
& \leq \|\boldtheta^n\|_{L^2(\Gamma_h)}\int_{t_{n-1}}^{t_n} \norm{\dot{\boldu}}_{L^2(\Gamma)} \leq \tau \norm{\dot{\boldu}}_{L^\infty([0,T], L^2(\Gamma))}\|\boldtheta^n\|_{L^2(\Gamma_h)}\\
& = C(\boldu)\tau \|\boldtheta^n\|_{L^2(\Gamma_h)}.
\end{split}
\end{equation}
From the Cauchy-Schwarz inequality and the estimate \eqref{dotensuresthatFD} for $\dot{\boldrho}$ we have
\begin{equation}
\label{p6FD}
\begin{split}
&\semin{\intgammah \bar{\partial}\boldrho^{n}:\boldtheta^{n}} \leq C\norm{\bar{\partial}\boldrho^{n}}_{L^2(\Gamma_h)}\norm{\boldtheta^{n}}_{L^2(\Gamma_h)} \\
&= \frac{C}{\tau} \norm{ \int_{t_{n-1}}^{t_n}\dot{\boldrho}}_{L^2(\Gamma_h)}\norm{\boldtheta^{n}}_{L^2(\Gamma_h)}\leq \frac{C}{\tau}\norm{\boldtheta^{n}}_{L^2(\Gamma_h)}\int_{t_{n-1}}^{t_n}\norm{\dot{\boldrho}}_{L^2(\Gamma_h)}\\
& \leq C\norm{\dot{\boldrho}}_{L^\infty([0,T], L^2(\Gamma_h))}\norm{\boldtheta^{n}}_{L^2(\Gamma_h)} \leq  C(\boldu)h^2\norm{\boldtheta^{n}}_{L^2(\Gamma_h)}.
\end{split}
\end{equation}
From the estimate \eqref{quadratureerror} for $\varepsilon_h$, the estimate \eqref{dotensuresthatFD} for $\boldrho$, the equivalences \eqref{equivalence1}, \eqref{equivalence2} and \eqref{equivalence4}, we obtain
\begin{equation}
\label{p7FD}
\begin{split}
&\semin{\varepsilon_h(\bar{\partial}\boldrho^{n},\boldtheta^{n})} \leq Ch^2\|\bar{\partial}\boldrho^{n}\|_{H^2_h(\Gamma_h)}\|\boldtheta^{n}\|_{H^1(\Gamma_h)}\\
&\leq \frac{Ch^2}{\tau}\norm{\boldtheta^{n}}_{H^1(\Gamma_h)}\int_{t_{n-1}}^{t_n}\|\dot{\boldrho}\|_{H^2_h(\Gamma_h)} \\
&\leq Ch^2\|\dot{\boldrho}\|_{L^\infty([0,T],H^2_h(\Gamma))}\norm{\boldtheta^{n}}_{H^1(\Gamma_h)}\\
& = Ch^2 (\|\dot{\boldrho}\|_{L^\infty([0,T], H^1(\Gamma_h))} + |\dot{\boldrho}|_{L^\infty([0,T],H^2_h(\Gamma_h))})\norm{\boldtheta^{n}}_{H^1(\Gamma_h)}\\
& = Ch^2(\|\dot{\boldrho}\|_{L^\infty([0,T],H^1(\Gamma_h))} + |\dot{\boldu}^{-\ell}|_{L^\infty([0,T], H^2_h(\Gamma_h))})\norm{\boldtheta^{n}}_{H^1(\Gamma_h)}\\
& \leq Ch^2(C(\boldu)h + (1+h)\norm{\dot{\boldu}}_{L^\infty([0,T], H^2(\Gamma))})\norm{\boldtheta^{n}}_{H^1(\Gamma_h)}\\
& \leq C(\boldu)(h^2+h^3)\norm{\boldtheta^{n}}_{H^1(\Gamma_h)}.
\end{split}
\end{equation}
From the Cauchy-Schwarz inequality and \eqref{equivalence1} we have
\begin{equation}
\label{p8FD}
\begin{split}
&\left|\intgammah (\bar{\partial}-\partial_t)\boldu^{-\ell,n}:\boldtheta^{n}\right| \leq C\norm{(\bar{\partial}-\partial_t)\boldu^{n}}_{L^2(\Gamma)}\norm{\boldtheta^{n}}_{L^2(\Gamma_h)}\\
& \leq \frac{C}{\tau}\norm{\boldtheta^{n}}_{L^2(\Gamma_h)}\int_{t_{n-1}}^{t_n}\norm{\dot{\boldu}(t)-\dot{\boldu}(t_n)}_{L^2(\Gamma)}\mathrm{d}t\\
& \leq \frac{C}{\tau}\norm{\boldtheta^{n}}_{L^2(\Gamma_h)}\int_{t_{n-1}}^{t_n}\int_t^{t_n}\norm{\ddot{\boldu}(s)} \mathrm{d}s\mathrm{d}t\\
& \leq C\tau \norm{\ddot{\boldu}}_{L^\infty([0,T],L^2(\Gamma))}\norm{\boldtheta^{n}}_{L^2(\Gamma_h)} = C(\boldu)\tau\norm{\boldtheta^{n}}_{L^2(\Gamma_h)},
\end{split}
\end{equation}
where we have exploited the assumption that $\ddot{\boldu}\in L^{\infty}([0,T],L^2(\Gamma))$. The Cauchy-Schwarz inequality, \eqref{equivalence1}, the geometric estimate \eqref{geometricalestimatedelta} and the stability bound \eqref{firststability} yield
\begin{equation}
\label{p9FD}
\begin{split}
&\semin{\intgamma \left( 1-\frac{1}{\delta_h^\ell}\right) \boldu^n:\boldtheta^{\ell,n}} \leq \norm{1-\frac{1}{\delta_h^\ell}}_{L^\infty(\Gamma)}\norm{\boldu^n}_{L^2(\Gamma)}\norm{\boldtheta^n}_{L^2(\Gamma_h)}\\
&\leq Ch^2\norm{\boldtheta^n}_{L^2(\Gamma_h)}.
\end{split}
\end{equation}
Combining \eqref{p1fully}-\eqref{p9FD}, using \eqref{equivalence} and Young's inequality we get
\begin{equation}
\label{p10FD}
\begin{split}
&\frac{1}{2\tau}( \|\boldtheta^{n}\|_h^2 - \|\boldtheta^{n-1}\|_h^2) + m\|\nablagammah \boldtheta^{n}\|_{L^2(\Gamma_h)}\\
&\leq C(\boldu)(h^2+h^3+\tau+\|\boldtheta^{n-1}\|_h)\norm{\boldtheta^{n}}_{H^1(\Gamma_h)}\\
&\leq C(\boldu,m)(h^4+h^6+\tau^2+\|\boldtheta^{n-1}\|_h^2) + m\norm{\boldtheta^{n}}_{H^1(\Gamma_h)}^2,
\end{split}
\end{equation}
where $m=\min_{k=1,\dots,r}\{d_k\}$, from which, canceling $\|\nablagammah\boldtheta^{-\ell,n}\|_{L^2(\Gamma_h)}$ on both sides of \eqref{p10FD}, and using \eqref{equivalence}, we have that
\begin{equation}
\label{almostfinishedFD}
\norm{\boldtheta^{n}}_h^2 \leq (1+C(\boldu)\tau)\|\boldtheta^{n-1}\|_h^2 + C(\boldu)\tau (h^4+h^6+\tau^2).
\end{equation}
By repeatedly applying \eqref{almostfinishedFD}, taking into account the assumption $\|\boldtheta^0\|_{L^2(\Gamma)}\leq Ch^2$, and then using \eqref{equivalence}, \eqref{equivalence1}, we obtain
\begin{equation*}
\|\boldtheta^{\ell,n}\|_{L^2(\Gamma)}^2 \leq C(\boldu)(h^4+h^6+\tau^2),
\end{equation*}
which yields the desired result.
\end{proof}
\end{theorem}
\noindent
The previous theorems imply that our semi- and fully-discrete methods exhibit optimal convergence rates, that is to say quadratic in the mesh size and linear in the time step.

\section{Numerical tests}
\label{sec:numericaltest}
In this section we solve some test problems numerically to show that the LSFEM combined with the IMEX Euler in time:
\begin{itemize}
\item exhibits the optimal convergence rate predicted in Theorem \ref{thm: FDerroranalysis} (Experiments \ref{sec:heatconvergence}, \ref{sec:exp4});
\item fulfills the discrete maximum principle for the homogeneous heat equation, whilst the standard SFEM does not (Experiment \ref{sec:exp2});
\item preserves the invariant rectangles of reaction-diffusion systems, whilst the standard SFEM does not (Experiment \ref{sec:exp3}).
\end{itemize}
The simulations have been carried out using MATLAB. The linear system arising at each timestep is solved with MATLAB's "backslash" command. The code is available on request.

\subsection{{\it Experiment 1}: The linear heat equation and its convergence}
\label{sec:heatconvergence}
\begin{figure}[h]
\centering{
\includegraphics[scale=0.65]{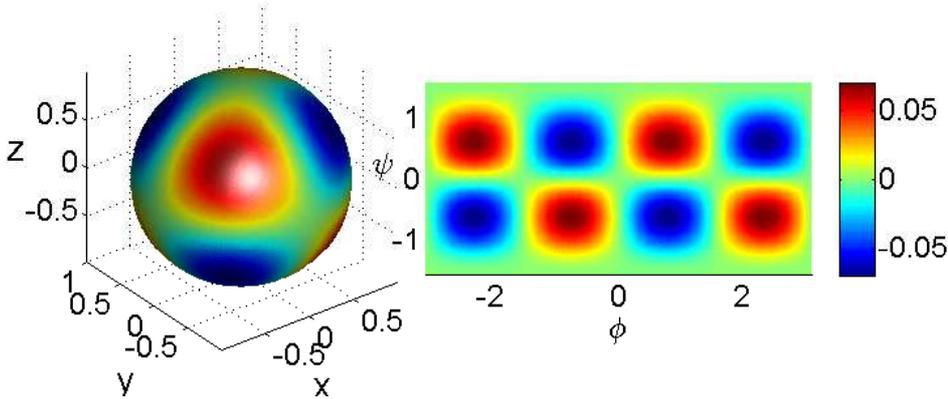}}
\caption{{\it Experiment 1}: The LSFEM solution corresponding to the linear heat equation \eqref{experiment1} with $d = \frac{1}{24}$ and $\beta = \frac{1}{2}$ obtained on a Delaunay mesh with $N=16962$ nodes and time step $\tau = 1.6$e-3 at $T=1$ (left) and its planar projection through spherical coordinates (right).}\label{fig:experiment1}
\end{figure}

In this experiment we solve the parabolic equation \eqref{systemunderstudy} in the linear case $\alpha =1$ on the unit sphere $\Gamma = \{(x,y,z)\in\mathbb{R}^3 | x^2+y^2+z^2=1\}$:
\begin{equation}
\label{experiment1}
\begin{cases}
&\dot{u} - d\Delta_\Gamma u = -\beta u,\\
&u_0(x,y,z) = xyz,\qquad (x,y,z)\in \Gamma,
\end{cases}
\end{equation}
with $d = \frac{1}{24}$ and $\beta = \frac{1}{2}$, to test the convergence rate of the LSFEM method. The exact solution of \eqref{experiment1} is
\begin{equation*}
u(x,y,z,t) = xyze^{-t},\qquad (x,y,z)\in\Gamma,\ t\geq 0.
\end{equation*}
In this experiment, as well as in the following ones, the problem is solved on a sequence of eight meshes $\Gamma_i$, $i=0,\dots,7$ with decreasing meshsizes $h_i\approx \sqrt{2}^{-i} h_0$ and corresponding time steps $\tau_i = 2^{-i} \tau_0$ (see parameter values in Tab. \ref{tab:experiment1}), so that $\tau_i$ is approximately proportional to $h_i^2$ in order to reveal the quadratic convergence, with respect to the mesh size, of the method. All of the $\tau_i$ fulfill the stability condition given in Theorem \ref{thm:fdmaximumprinciple}. For every $i=0,\dots,7$ the $L^\infty([0,T],L^2(\Gamma_h))$ error between the numerical solution $U$ and the interpolant $I_h(u)$ of the exact solution is measured and the numerical results are reported in Table \ref{tab:experiment1}. The lumped solution at the final time $T=1$ obtained on the finest mesh is shown in Figure \ref{fig:experiment1} (left), as well as its planar projection through spherical coordinates
\begin{equation*}
x = \cos\phi\cos\psi,\quad y = \cos\phi\sin\psi, \quad z = \sin\psi,\qquad (\phi,\psi)\in [-\pi,\pi]\times \left[-\frac{\pi}{2}, \frac{\pi}{2}\right].
\end{equation*}
In this test example we observe that the lumped SFEM is more accurate than the standard SFEM and the predicted second order convergence in space is attained.

\begin{table}
\caption{{\it Experiment 1}: Comparisons of the convergence analysis between the SFEM and the LSFEM for the linear heat equation \eqref{experiment1} with $d = \frac{1}{24}$ and $\beta = \frac{1}{2}$.}\label{tab:experiment1}
\input{table1.tex}
\end{table}

\subsection{{\it Experiment 2}: The homogeneous heat equation and the maximum principle}
\label{sec:exp2}
\begin{figure}[h]
\centering{
\includegraphics[scale=0.65]{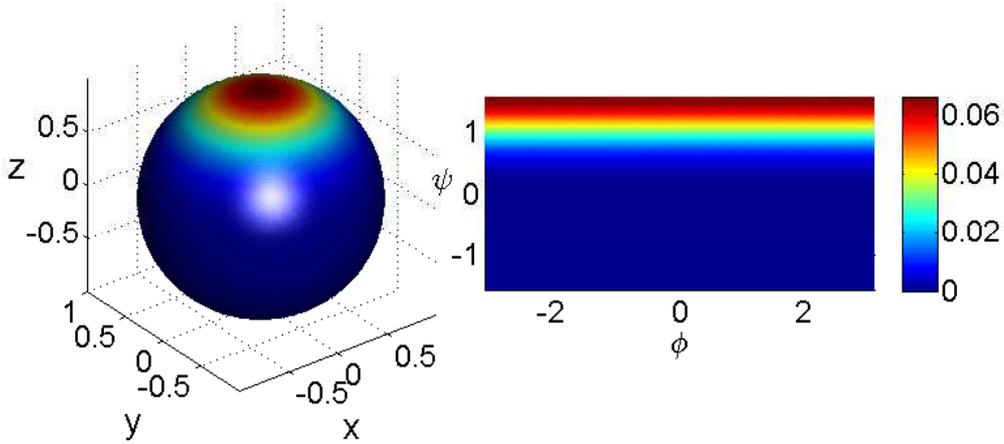}}
\caption{{\it Experiment 2}: The LSFEM solution corresponding to the homogeneous heat equation \eqref{systemunderstudy} with $\beta=0$ and initial datum \eqref{initialheatmaximumprinciple} obtained on a Delaunay mesh with $N=16962$ nodes and time step $\tau = 1.6e$-3 at $T=1$ (left) and its planar representation through spherical coordinates (right).}\label{fig:experiment2}
\end{figure}

We solve the parabolic equation \eqref{systemunderstudy} for the homogeneous case $\beta=0$ on the unit sphere $\Gamma$ with $d=0.1$ and the nonnegative compactly supported $H^1(\Gamma)$ initial datum
\begin{equation}
\label{initialheatmaximumprinciple}
u_0(x,y,z) = 
\begin{cases}
&\sqrt{1-\frac{x^2+y^2}{0.04}}\qquad \text{if}\ x^2+y^2 \leq 0.04,\ z>0,\\
&\hspace{1cm} 0\hspace{0.8cm} \qquad \text{elsewhere},
\end{cases}
\end{equation}
The minima of the computed numerical solution, obtained for every choice of $(h,\tau)$, are reported in Table \ref{tab:experiment2}. In Figure \ref{fig:experiment2} we show the LSFEM solution obtained on the finest mesh at the final time. This experiment confirms our findings, as the LSFEM fulfills the discrete maximum principle, while the standard SFEM violates the maximum principle as illustrated in Table \ref{tab:experiment2}.

\begin{table}
\caption{{\it Experiment 2}: Discrete maximum principle analysis: comparisons between SFEM and LSFEM for the homogeneous heat equation \eqref{systemunderstudy} with $\beta=0$ and initial datum \eqref{initialheatmaximumprinciple}.}\label{tab:experiment2}
\input{table2.tex}
\end{table}

\subsection{{\it Experiment 3}: Reaction-diffusion system and the preservation of the invariant rectangle}
\label{sec:exp3}
In this experiment we consider the  reaction-diffusion system with Rosenzweig-MacArthur kinetics \cite{gonzalez2003dynamic, garvie2007} 
\begin{equation}
\label{numericalrd}
\begin{cases}
u_t -d_1\Delta_{\Gamma} u = au(1-u) -b\frac{uv}{u+\alpha},\\
v_t -d_2\Delta_{\Gamma} v = c\frac{uv}{u+\alpha} -dv,
\end{cases}
\end{equation}
on the unit sphere $\Gamma$, where $\alpha, a,b,c,d$ are positive constants.

This system has been numerically solved in \cite{garvie2007} on a planar domain with LFEM in combination with an implicit Euler time discretization. However, since the theory developed in \cite{garvie2007} addresses a problem on domains of more general dimension ($n\leq 3$) there is no discrete maximum principle and the authors consider modified kinetics to ensure the positiveness of the numerical solution. The present example shows that, on two dimensional manifolds, lumping guarantees the preservation of the invariant region without needing modified kinetics.

\begin{figure}
\centering{
\includegraphics[scale=0.35]{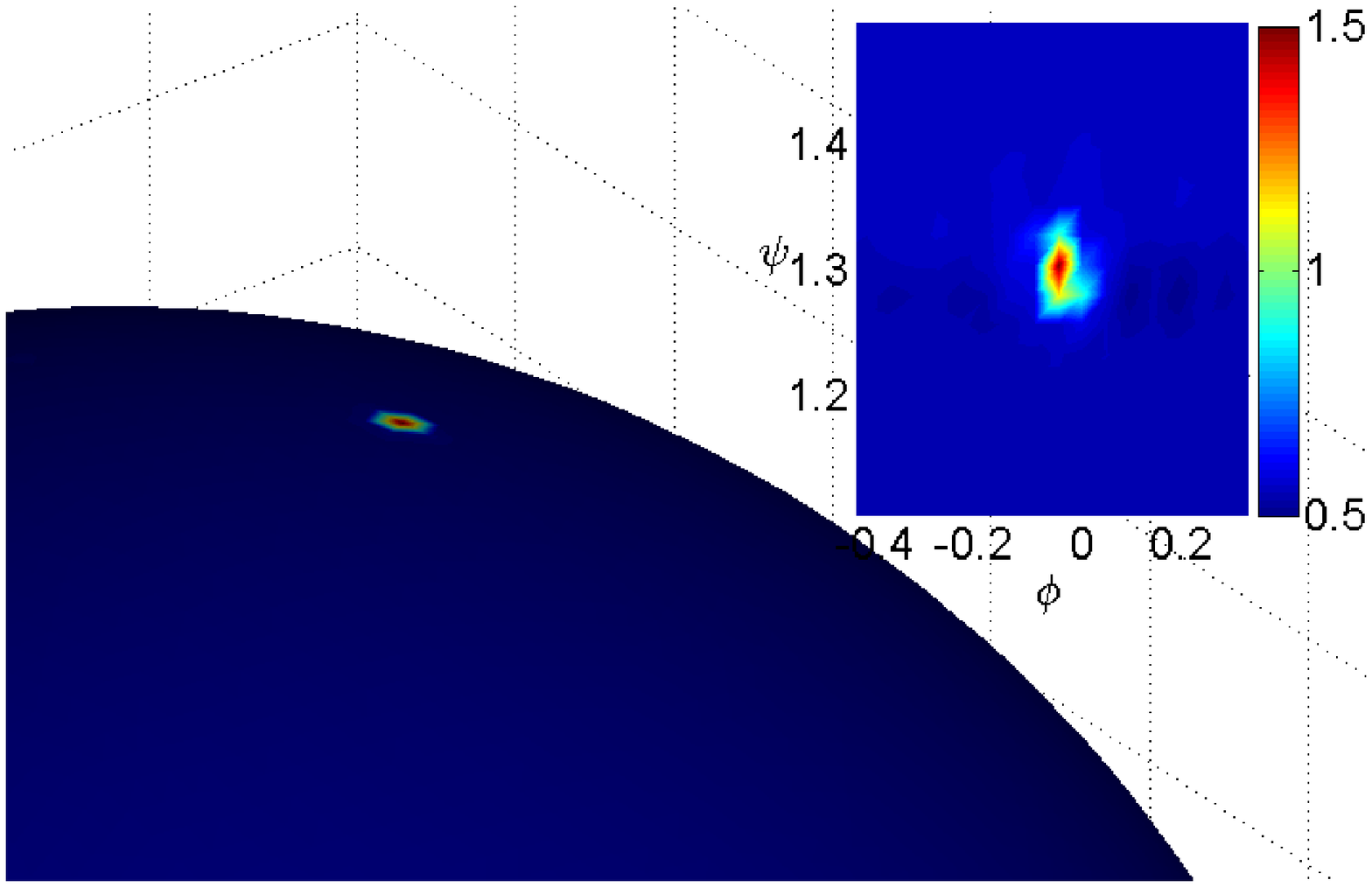}
\includegraphics[scale=0.35]{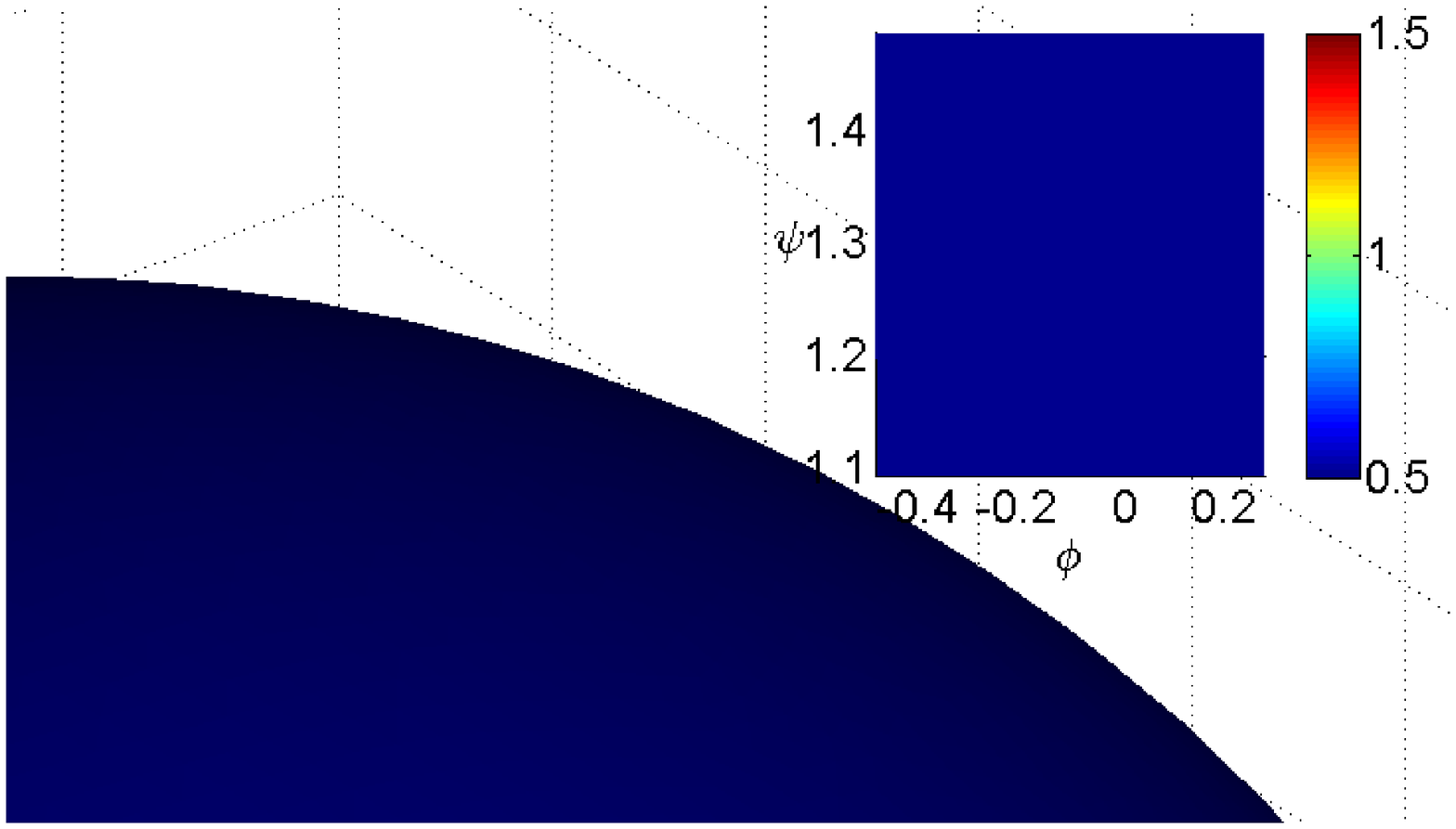}
}\caption{{\it Experiment 3}: Component $v$ of the numerical solution of \eqref{numericalrd} with $d_1= d_2 = 1\text{e-2}$, $\alpha=1\text{e-3}$, $a=10$, $b=1\text{e-2}$, $c=d=1$ and initial datum \eqref{exp3initialu}-\eqref{exp3initialv} with $r=0.2, \varepsilon = 1$e-7 obtained on a mesh with $N=16962$ gridpoints at $\bar{\bar{t}}=0.043$ by SFEM (left panel) and LSFEM (right panel). The zoom insets show that the solution of the SFEM method overcomes the threshold 0.5 of the invariant region, whilst the LSFEM does not.}\label{fig:experiment31}
\end{figure}

When $c=d$ and $0<\alpha<\frac{1}{\sqrt{2}}$ for every $0<\varepsilon<1-2\sqrt{a}$, the rectangle
\begin{equation}
\label{numericalinvariantregion}
\Sigma := \left[\varepsilon, 1\middle] \times \middle[0, \frac{a\alpha}{2b}\right]
\end{equation}
is an invariant region for \eqref{numericalrd}, see for instance the analysis in \cite{gonzalez2003dynamic}. An easy way to see this is to observe that, for every $\varepsilon,\varepsilon'>0$, the rectangle 
\begin{equation*}
\Sigma_1 := \left[\varepsilon, 1+ \frac{\varepsilon' a \alpha}{b}\middle] \times \middle[-\varepsilon', \frac{a\alpha}{2b}\right]
\end{equation*}
fulfills condition \eqref{smollercondition}. Then, since the intersection of invariant regions is still invariant, also $\Sigma$ is invariant  for \eqref{numericalrd}. The $H^1(\Gamma)$ initial datum
\begin{align}
\label{exp3initialu}
&u_0(x,y,z) = \begin{cases}
&\varepsilon + (1-\varepsilon)\sqrt{1-\frac{x^2+y^2}{r^2}}\qquad \text{if}\ x^2+y^2 \leq r^2,\ z>0,\\
&\varepsilon\qquad \text{elsewhere},
\end{cases}\\
\label{exp3initialv}
&v_0(x,y,z) = \frac{a\alpha}{2b}, \qquad \forall\ (x,y,z)\in\Gamma,
\end{align}
with $0<r<1$, is contained in the invariant region $\Sigma$. Furthermore, for $0<\alpha<1$, it is easy to verify that, on $\Sigma$, the Lipschitz constants $L_1$ and $L_2$ of the kinetics in \eqref{numericalrd} fulfill
\begin{equation*}
L_1 < \sqrt{2}\left(3a+\frac{b}{2\alpha}\right),\quad \text{and}  \quad L_2 < \sqrt{2}\left(\frac{c}{2\alpha} + \frac{d}{2}\right).
\end{equation*}
In the following we choose $d_1= d_2 = 1\text{e-2}$, $\alpha=1\text{e-3}$, $a=10$, $b=1\text{e-2}$, $c=d=1$, $r=0.2$, and $\varepsilon = 1\text{e-7}$. With these settings the invariant region \eqref{numericalinvariantregion} becomes
\begin{equation}
\label{particularinvariantregion}
\Sigma = [\text{1e-7},1]\times \left[0,\frac{1}{2}\right],
\end{equation}
and the stability condition \eqref{fdstabilitycondition} on the time step is fulfilled if we choose
\begin{equation}
\tau \leq \bar{\tau} := \frac{1}{\sqrt{2}\max\left\{\left(3a+\frac{b}{2\alpha}\right), \left(\frac{c}{2\alpha} + \frac{d}{2}\right)\right\} } = 1.4e\text{-3}.
\end{equation}
We thus solve the problem on the same sequence of spatial meshes considered in the previous experiments, with a fixed time step $\bar{\tau}= 1e$-3 and final time $T=5$. In Tables \ref{tab:experiment31}-\ref{tab:experiment32} we show the minima and the maxima of the components of the computed numerical solution: we observe that the LSFEM solution preserves $\Sigma$, whilst the SFEM one violates $\Sigma$ on all of the considered meshes. Furthermore, the SFEM exhibits a stability threshold: the numerical solution blows up on meshes $i=0,\dots,4$, while it stays bounded on the finer meshes $i=5,6,7$. In fact, on the latter meshes, the absolute minima and maxima are attained within $\bar{t} := 0.121$, while the final computational time is $T=5$. In Figure \ref{fig:experiment31} we show the $v$ component of both the SFEM and LSFEM solutions, computed on mesh $i=7$, at the time $\bar{\bar{t}} := 0.043$ in which the $v$ component of the SFEM solution attains its absolute minimum. It is evident that the solution of the SFEM method overcomes the threshold 0.5 of the invariant region, whilst the LSFEM does not.
\begin{table}[h]
\caption{{\it Experiment 3}: Invariance analysis for the SFEM solution of \eqref{numericalrd} with parameters and initial datum as stated in Fig. \ref{fig:experiment31}. The solution clearly blows up on the first five meshes. On the three finest meshes the numerical solution stays bounded, though still violating the invariant region \eqref{particularinvariantregion}.}\label{tab:experiment31}
\input{table3.tex}
\end{table}
\begin{table}
\caption{{\it Experiment 3}: Invariance analysis for the LSFEM solution of \eqref{numericalrd} with parameters and initial datum as stated in Fig. \ref{fig:experiment31}. The solution stays in the invariant rectangle $[\text{1e-7},1]\times\left[0,\frac{1}{2}\right]$ on all of the considered meshes.}\label{tab:experiment32}
\input{table4.tex}
\end{table}

\subsection{{\it Experiment 4}: Reaction-diffusion system with {\it activator-depleted} kinetics and its convergence}
\label{sec:exp4}
\begin{figure}
\centering{
\includegraphics[scale=0.65]{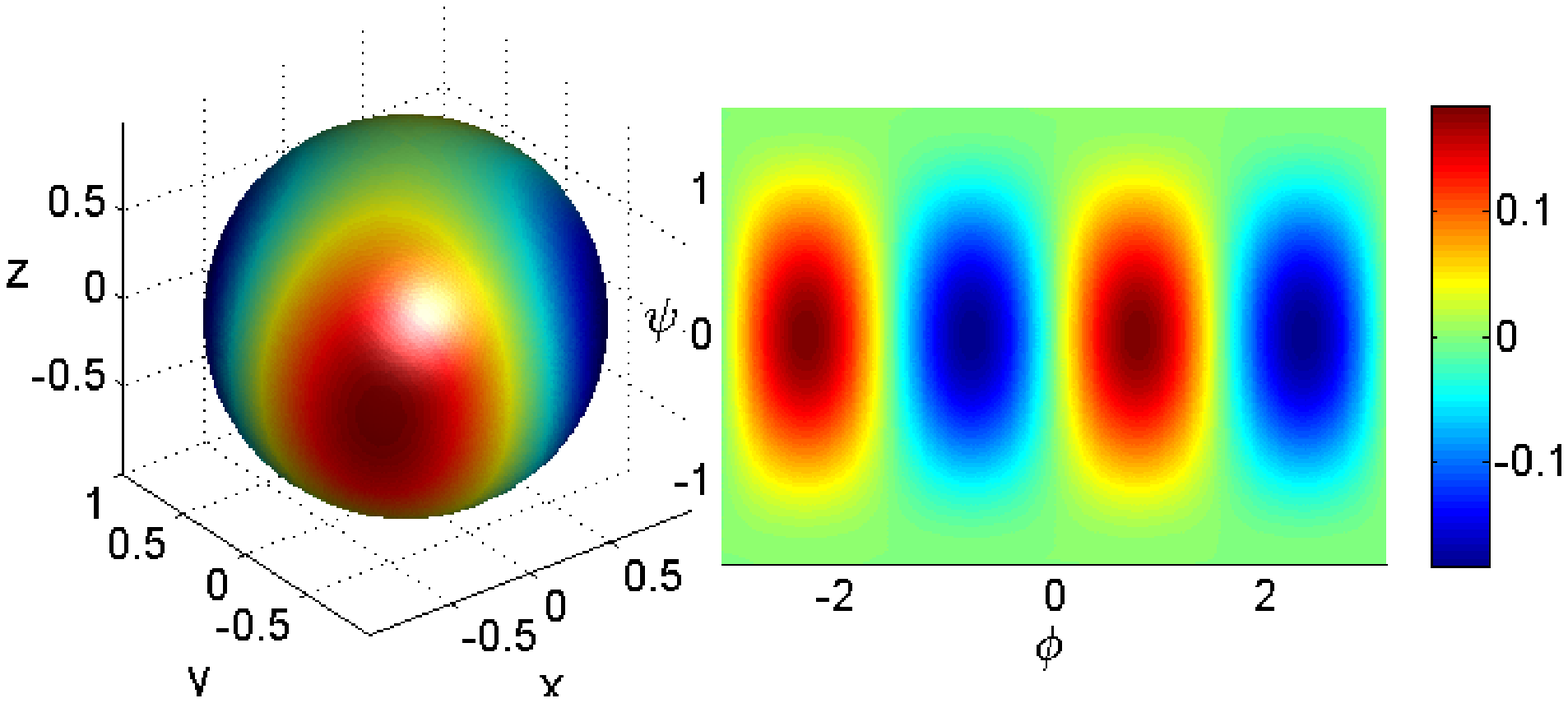}
\includegraphics[scale=0.45]{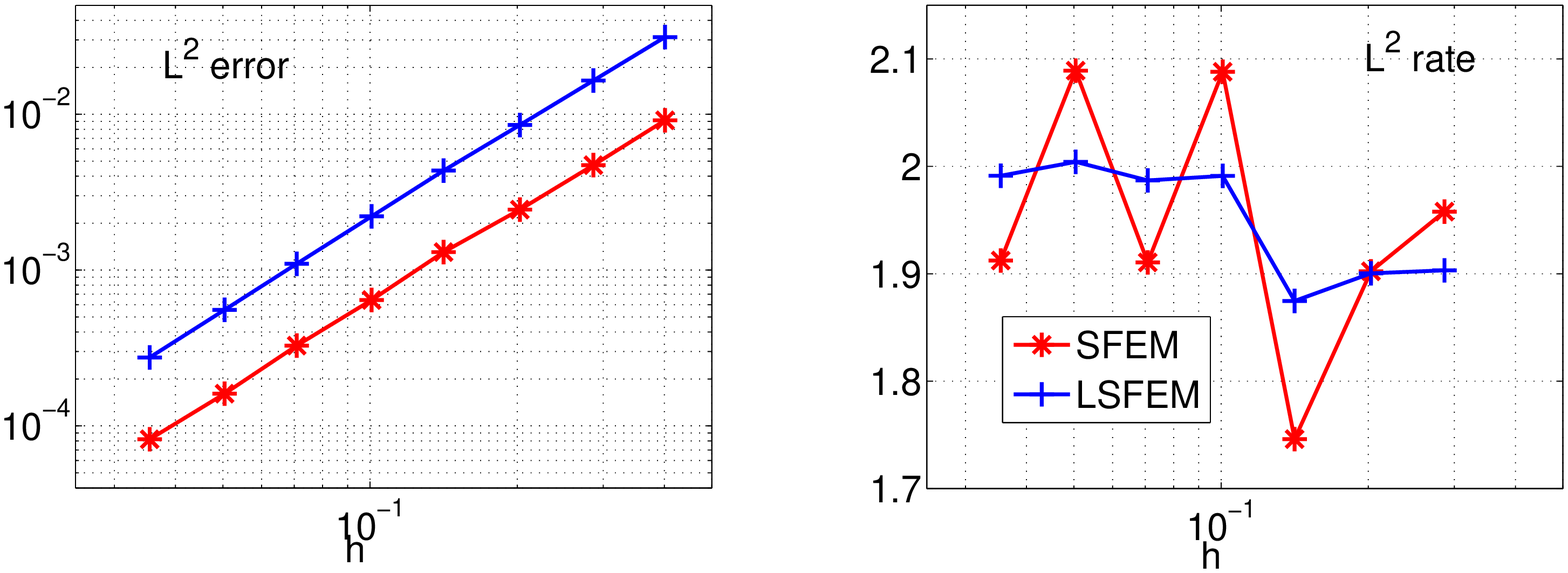}}
\caption{Top row: the $u$-component of the LSFEM solution corresponding to the reaction-diffusion system with {\it activator-depleted} substrate kinetics \eqref{numericalrdconvergence1}-\eqref{exp4forcingterms} with $a=b=1$, $d_1=\frac{1}{6}$, $d_2=\frac{1}{12}$ and initial condition \eqref{exp4initialcond} obtained on a mesh with $N=16962$ nodes and time step $\tau = 1.6e$-3 at $T=1$ and their corresponding planar projections through spherical coordinates. Bottom row: convergence analysis of the SFEM and LSFEM. As predicted, the LSFEM retains the quadratic convergence rate of the SFEM.}\label{fig:experiment41}
\end{figure}

In this example, we test the convergence rate of the method on a reaction-diffusion system on the unit sphere $\Gamma$ with well-studied {\it activator-depleted} substrate kinetics\cite{gierer1972theory,prigogine1968symmetry,schnakenberg1979simple,murray2001} with an additional forcing term on the right hand side:
\begin{equation}
\label{numericalrdconvergence1}
\begin{cases}
u_t -d_1\Delta_{\Gamma} u = a - u + u^2v + f_1(x,y,z,t),\\
v_t -d_2\Delta_{\Gamma} v = b - u^2v + f_2(x,y,z,t),
\end{cases}
\end{equation}
with the functions $f_1(\boldx,t), f_2(\boldx,t)$ chosen in such a way that the exact solution is known at all times. Although this example is beyond the scope of the present work, due to the space
and time dependence of the reaction terms, we include it merely as a numerical test.\\
We choose $a=b=1$, $d_1=\frac{1}{6}$, $d_2=\frac{1}{12}$, 
\begin{equation}
\label{exp4forcingterms}
\begin{cases}
&f_1(x,y,z,t) = xye^{-t}(1 + x^2y^2e^{-2t}) -a;\\
&f_2(x,y,z,t) = -x^3y^3ze^{-t}-b,
\end{cases}
\end{equation}
and the following initial condition:
\begin{equation}
\label{exp4initialcond}
\begin{cases}
&u_0(x,y,z) = xy,\\
&v_0(x,y,z) = -xyz,
\end{cases}
\qquad \forall\ (x,y,z)\in\Gamma.
\end{equation}
In this case, the exact solution is given by
\begin{equation*}
\begin{cases}
&u(x,y,z,t) = xye^{-t},\\
&v(x,y,z,t) = -xyze^{-t},
\end{cases}
\qquad \forall\ (x,y,z)\in \Gamma,\ \forall\ t\geq 0.
\end{equation*}
We solve the problem on the same sequence of meshes and time steps considered in {\it Experiment 1}, with final time $T=1$, for both the SFEM and the LSFEM, where the contributions due to the forcing terms $f_i$, $i=1,2$ are approximated with the standard and the lumped quadrature rule given by
\begin{align*}
\intgammah I_h(f_i)\chi_j, \quad \text{and} \quad 
\intgammah I_h(f_i\chi_j), \quad \forall\ i=1,\dots, N,
\end{align*}
respectively. We observe that the standard quadrature rule is exact for piecewise linear functions, whilst the lumped one is only exact when the product of the functions is piecewise linear. For this reason, the LSFEM is expected to produce larger errors than the SFEM. The $L^2$ errors and experimental convergence rates are plotted in Fig. \ref{fig:experiment41} together with the LSFEM solution obtained on the finest mesh at the final time $T=1$. As expected, the LSFEM exhibits slightly larger errors than the SFEM. Nonetheless, they have the same convergence rate, in agreement with our  theoretical findings.

\section{Conclusions}
In Section \ref{sec:heatproblem} we have considered a lumped surface finite element method (LSFEM), for a class of semilinear parabolic problems on surfaces, by extending its planar counterpart \cite{thomee1980} inspired by the ideas in \cite{acta}. Time discretisation is carried out by applying the IMEX Euler method in time. We have shown in Theorem \ref{thm:sdmaximumprinciple} that the spatially discrete problem fulfills a discrete maximum principle. In particular, we have proved that: no restriction on the timestep is required in the homogeneous case (thus extending the result of \cite{thomeebook} to surfaces); the time step restriction \eqref{heattimesteprestriction} is required in the presence of the nonlinear reaction terms in \eqref{systemunderstudy}.
In Section \ref{sec:rdsysyems} we have applied the LSFEM to general systems of arbitrarily many reaction-diffusion equations. In analogy to the continuous setting (see \cite{smollerarticle}), in Theorem \ref{thm:sdinvariantregion} we have shown that, under the sole assumption of Delaunay regularity for the mesh, the strictly inward flux condition \eqref{smollercondition} is sufficient for a rectangle in the phase space to be invariant for the spatially discrete scheme. For the fully-discrete problem arising from IMEX Euler we have shown in Theorem \ref{thm:fdinvariantregion} that, under the time step restriction \eqref{fdstabilitycondition} involving the Lipschitz constants of the reaction kinetics,  condition \eqref{smollercondition} is not only still sufficient to ensure a hyper-rectangle to be invariant, but can be even weakened by requiring non-outward fluxes \eqref{cond2}. To the best of the authors' knowledge, Theorems \ref{thm:sdinvariantregion} and \ref{thm:fdinvariantregion} are a novelty even on planar domains.\\
For both the semi- and fully-discrete formulations of the reaction-diffusion systems considered in Section \ref{sec:rdsysyems}, including the parabolic problem of Section \ref{sec:heatproblem} as a special case, an optimal $L^2(\Gamma)$ error bound has been proven in Section \ref{sec:erroranalysis}. The numerical examples in Section \ref{sec:numericaltest} confirm our theoretical findings. The usefulness of LSFEM is illustrated in Experiments \ref{sec:exp2} and \ref{sec:exp3}. In particular, we have shown that in the absence of lumping the numerical solutions of the homogeneous heat equation violates the maximum principle (Exp 5.2) and the numerical solution of a classical predator-prey model blows-up instead of being bounded by the invariant rectangle.\\
Emerging applications encourage the extension of the present study to the case of  \emph{evolving} surfaces, which is beyond the scope of this work and will be addressed in future studies.

\section*{Acknowledgements}
This work (AM, CV) is partly supported by the EPSRC grant number EP/J016780/1 and the
Leverhulme Trust Research Project Grant (RPG-2014-149). The authors (MF, AM, IS CV) would like to thank the Isaac Newton Institute for
Mathematical Sciences for its hospitality during the programme [Coupling Geometric PDEs with Physics for Cell Morphology, Motility and Pattern
Formation] supported by EPSRC Grant Number EP/K032208/1.  AM acknowledges funding from the European Union Horizon 2020 research and
innovation programme under the Marie Sk\l{}odowska-Curie grant agreement No 642866 and
was partially supported by a grant from the Simons Foundation.

\input{bibliography.bbl}
\end{document}

%% file: figure1.tex
\begin{tikzpicture}
\draw[color =red] (0,0)--(1,3);
\draw (1,3)--(3,-1)--(0,0)--(-2,-1)--(1,3);
\draw[color = red, domain = 117:160] plot ({3+cos(\x)/2}, {-1+sin(\x)/2});
\draw[color = red, domain = 117:160] plot ({3+cos(\x)/1.8}, {-1+sin(\x)/1.8});
\draw[color = red, domain = 25:54] plot ({-2+cos(\x)/2}, {-1+sin(\x)/2});
\draw[color = red, domain = 25:54] plot ({-2+cos(\x)/1.8}, {-1+sin(\x)/1.8});
\draw[dotted] (3,2.5)--(1,0.5);
\draw[dotted] (-2.5,2)--(-0.3,0.5);
\node[label = right:{$K_2$}] at (3,2.5) {};
\node[label = left:{$K_1$}] at (-2.5,2) {};
\node[label = right:{$\red{\alpha_2}$}] at (3,-1) {};
\node[label = left:{$\red{\alpha_1}$}] at (-2,-1) {};
\node[label = right:{$\red{e}$}] at (0.5,1.5) {};
\end{tikzpicture}

%% file: table1.tex
\begin{tabular}{c c c c c c c}
\hline\noalign{\smallskip}
\multicolumn{3}{c}{ } & \multicolumn{2}{c}{SFEM} & \multicolumn{2}{c}{LSFEM}\\
\noalign{\smallskip}\hline\noalign{\smallskip}
$i$ & $N$ & $h$ & $L^2$ error & $L^2$ rate & $L^2$ error & $L^2$ rate \\
\noalign{\smallskip}\hline\noalign{\smallskip}
0 & 126 & 4.013e-01 & 6.100e-03 & -         & 3.061e-03 & -        \\
1 & 258 & 2.863e-01 & 3.129e-03 & 1.977 & 1.846e-03 & 1.498 \\
2 & 516 & 2.026e-01 & 1.594e-03 & 1.951 & 1.095e-03 & 1.510 \\
3 & 1062 & 1.414e-01 & 7.899e-04 & 1.953 & 5.444e-04 & 1.945 \\
4 & 2094 & 1.007e-01 & 3.966e-04 & 2.030 & 3.025e-04 & 1.731 \\
5 & 4242 & 7.082e-02 & 2.013e-04 & 1.925 & 1.401e-04 & 2.184 \\
6 & 8370 & 5.041e-02 & 1.003e-04 & 2.049 & 7.671e-05 & 1.773 \\
7 & 16962 & 3.542e-02 & 5.063e-05 & 1.938 & 3.529e-05 & 2.200 \\
\noalign{\smallskip}\hline
\end{tabular}

%% file: table2.tex
\begin{tabular}{c c c c c}
\hline\noalign{\smallskip}
$i$ & $N$ & $h$ & $\min_{\Gamma_h\times [\tau,1]} U$ SFEM& $\min_{\Gamma_h\times [\tau,1]} U$ LSFEM\\
\noalign{\smallskip}\hline\noalign{\smallskip}
0 &126 & 4.013e-01&-3.454e-04& 7.016e-09\\
1 & 258 & 2.863e-01&-4.695e-06& 4.812e-12\\
2 & 516 & 2.026e-01&-1.299e-03& 1.213e-16\\
3 & 1062 & 1.414e-01&-2.123e-07& 2.746e-23\\
4 & 2094 & 1.007e-01&-7.546e-04& 3.142e-32\\
5 & 4242 & 7.082e-02&-1.037e-05& 1.816e-45\\
6 & 8370 & 5.041e-02&-4.163e-04& 5.324e-64\\
7 & 16962 & 3.542e-02&-1.254e-04& 3.126e-90\\
\noalign{\smallskip}\hline
\end{tabular}

%% file: table3.tex
\begin{tabular}{c c c c c c c}
\hline\noalign{\smallskip}
$i$ & $N$ & $h$&$\min_{\Gamma_h\times [\tau,5]}$ U&$\max_{\Gamma_h\times [\tau,5]}$ U&$\min_{\Gamma_h\times [\tau,5]}$ V&$\max_{\Gamma_h\times [\tau,5]}$ V\\
\noalign{\smallskip}\hline\noalign{\smallskip}
0 & 126 & 4.013e-01 & -9.876e+268 & 9.477e+264 & -3.965e-02 & 5.503e-01\\
1 & 258 & 2.863e-01 & -8.694e+255 & 1.344e+252 & -5.517e-02 & 1.798e+00\\
2 & 516 & 2.026e-01 & -1.749e+255 & 6.160e+251 & -2.227e-01 & 2.081e+00\\
3 & 1062 & 1.414e-01 & -3.237e+274 & 2.577e+271 & -4.924e-01 & 4.711e+00\\
4 & 2094 & 1.007e-01 & -1.309e+267 & 1.553e+264 & -3.862e+01 & 8.042e-01\\
5 & 4242 & 7.082e-02 & -1.653e-02   & 9.999e-01    & -3.013e+00 & 9.192e-01\\
6 & 8370 & 5.041e-02 & -1.317e-02   & 9.999e-01   & -6.236e-01  & 1.706e+00\\
7 & 16962 & 3.542e-02 & -1.440e-02   & 9.999e-01  & 1.403e-01 & 1.865e+00\\
\noalign{\smallskip}\hline
\end{tabular}

%% file: table4.tex
\begin{tabular}{c c c c c c c}
\hline\noalign{\smallskip}
$i$ & $N$&$h$&$\min_{\Gamma_h\times [\tau,5]}$ U&$\max_{\Gamma_h\times [\tau,5]}$ U&$\min_{\Gamma_h\times [\tau,5]}$ V&$\max_{\Gamma_h\times [\tau,5]}$ V\\
\noalign{\smallskip}\hline\noalign{\smallskip}
0 & 126 & 4.013e-01 & 1.005e-07 & 9.999e-01  &  1.403e-01  &  4.999e-01\\
1 & 258 & 2.863e-01 & 1.005e-07 & 9.999e-01  & 1.403e-01  & 4.999e-01\\
2 & 516 & 2.026e-01 & 1.005e-07 & 9.999e-01  & 1.403e-01 & 4.999e-01\\
3 & 1062 & 1.414e-01 & 1.005e-07 & 9.999e-01  & 1.403e-01 & 4.999e-01\\
4 & 2094 & 1.007e-01 & 1.005e-07 & 9.999e-01  & 1.403e-01 & 4.999e-01\\
5 & 4242 & 7.082e-02 & 1.005e-07 & 9.999e-01  & 1.403e-01  & 4.999e-01\\
6 & 8370 & 5.041e-02 & 1.005e-07 & 9.999e-01  & 1.403e-01 & 4.999e-01\\
7 & 16962 & 3.542e-02 & 1.005e-07 & 9.999e-01  & 1.403e-01 & 4.999e-01\\
\noalign{\smallskip}\hline
\end{tabular}